\documentclass[12pt, a4paper]{amsart}
\usepackage{fullpage}
\pdfoutput=1
\usepackage{amsmath}

\usepackage{amssymb}
\usepackage{stmaryrd}%for \longmapsfrom
\usepackage{amsthm}

\usepackage{amscd}
\usepackage[all]{xy}
\usepackage{url}

\usepackage[colorlinks]{hyperref}
\hypersetup{%
pdftitle={Quantum Cluster Variables via Serre Polynomials},
pdfauthor={Fan QIN},
pdfkeywords={Quantum Cluster Algebra, Quiver Grassmannians, Serre polynomials, Positivity Conjecture},
bookmarksnumbered,
pdfstartview={FitH},
breaklinks=true,
urlcolor=blue,
citecolor=blue,
}%
\usepackage[all]{hypcap} %to fix fraction numbers showing problems in package {hyperref}
\newtheorem{Thm}{Theorem}[subsection]
\newtheorem{Lem}[Thm]{Lemma}
\newtheorem{Prop}[Thm]{Proposition}
\newtheorem{Cor}[Thm]{Corollary}

\newtheorem{Rem}[Thm]{Remark}
\newtheorem{Def}[Thm]{Definition}

\newtheorem*{Def*}{Definition}
\newtheorem*{Thm*}{Theorem}

\newcommand{\Z}{\mathbb{Z}}
\newcommand{\N}{\mathbb{N}}
\newcommand{\Q}{\mathbb{Q}}
\newcommand{\C}{\mathbb{C}}
\newcommand{\F}{\mathbb{F}}

\newcommand{\Fp}{{\mathbb F}_p}
\renewcommand{\P}{\mathbb{P}}

\newcommand{\ie}{{\em i.e.}\ }
\newcommand{\cf}{{\em cf.}\ }
\newcommand{\eg}{{\em e.g.}\ }

\newcommand{\st}{{such that}\ }

\renewcommand{\tilde}[1]{\widetilde{#1}}

\newcommand{\opname}[1]{\mathsf{#1}}%change definition !!

\renewcommand{\mod}{\opname{mod}}%redefined!

\newcommand{\per}{\opname{per}}

\newcommand{\Dfd}{\cD_{fd}}

\newcommand{\pr}{\opname{pr}}

\newcommand{\ind}{\opname{ind}}

\newcommand{\add}{\opname{add}}

\newcommand{\ra}{\rightarrow}
\newcommand{\xra}{\xrightarrow}

\newcommand{\embed}{\hookrightarrow}
\newcommand{\iso}{\stackrel{_\sim}{\rightarrow}}

\newcommand{\Gr}{\opname{Gr}}

\newcommand{\dimv}{\underline{\dim}\,}

\newcommand{\Ker}{\opname{Ker}}
\renewcommand{\Im}{\opname{Im}}%redefined!
\newcommand{\Cok}{\opname{Cok}}

\newcommand{\sign}{\opname{sign}}

\newcommand{\Hom}{\opname{Hom}}
\newcommand{\RHom}{\opname{RHom}}

\newcommand{\Ext}{\opname{Ext}}

\newcommand{\ten}{\otimes}
\newcommand{\lten}{\overset{\textbf{L}}{\ten}} %redefined

\renewcommand{\deg}{\opname{deg}}

\newcommand{\K}{\operatorname{\mathbb{K}}\nolimits}
\newcommand{\lat}{\opname{lat}}
\newcommand{\cc}{{\mathcal C}}
\newcommand{\cd}{{\mathcal D}}

\newcommand{\Rq}{{\Z[q^{\pm\frac{1}{2}}]}}
\newcommand{\Rk}{{\Z[|\sK|^{\pm\frac{1}{2}}]}}
\newcommand{\Kz}{{\opname{K}_0}}
\newcommand{\Hf}{{\frac{1}{2}}}

\newcommand{\Ra}{\longrightarrow}

\newcommand{\Rm}[1]{{\longmapsto}}
\newcommand{\Lm}[1]{{\longmapsfrom}}
\newcommand{\xRa}[1]{\stackrel{#1}{\longrightarrow}}

\newcommand{\aL}{{L}}%abstract Lattice
\newcommand{\sch}{\opname{Sch}}%category of scheme

\newcommand{\mW}{{W}}%%%% module W
\newcommand{\mV}{{V}}
\newcommand{\mU}{{U}}
\newcommand{\mI}{{I}}
\newcommand{\mP}{{P}}
\newcommand{\mM}{{M}}
\newcommand{\mN}{{L}}%% originally N
\newcommand{\mBp}{{E}}%% originally B_+
\newcommand{\mBn}{{E'}}%%originally B_-
\newcommand{\tmW}{{\tW}}
\newcommand{\tmV}{{\tV}}
\newcommand{\tmU}{{\tU}}
\newcommand{\tmM}{{\tM}}
\newcommand{\tmN}{{\tL}}%changed
\newcommand{\tmBp}{{\tilde{E}}}%changed
\newcommand{\tmBn}{{\tilde{E'}}}%changed
\newcommand{\tmI}{{\tI}}
\newcommand{\tmP}{{\tP}}

\newcommand{\mA}{{A}}
\newcommand{\mB}{{B}}
\newcommand{\mC}{{C}}
\newcommand{\mD}{{D}}
\newcommand{\mE}{{N}}%changed
\newcommand{\mF}{{O}}%changed
\newcommand{\mL}{{J}}%changed
\newcommand{\mK}{{K}}
\newcommand{\mX}{{X}}
\newcommand{\mY}{{Y}}

\newcommand{\da}{{a}}%%% dimension of A
\newcommand{\db}{{b}}
\newcommand{\dc}{{c}}
\newcommand{\dd}{{d}}
\newcommand{\de}{{n}}%changed
\newcommand{\df}{{o}}%changed
\newcommand{\dm}{{m}}
\newcommand{\dn}{{l}}%changed
\newcommand{\dl}{{j}}%changed
\newcommand{\dk}{{k}}
\newcommand{\dx}{{x}}
\newcommand{\dy}{{y}}

\newcommand{\mG}{{G}}
\newcommand{\mH}{{H}}
\newcommand{\dg}{{g}}
\renewcommand{\dh}{{h}}
\newcommand{\dv}{{v}}

\newcommand{\dw}{{w}}

\newcommand{\mTheta}{{\Im\theta}}
\newcommand{\mPhi}{{\Im\phi}}

\newcommand{\cA}{{\mathcal A}}

\newcommand{\cC}{{\mathcal C}}
\newcommand{\cD}{{\mathcal D}}

\newcommand{\cF}{{\mathcal F}}

\newcommand{\cP}{{\mathcal P}}

\newcommand{\cR}{{\mathcal R}}

\newcommand{\cT}{{\mathcal T}}
\newcommand{\cU}{{\mathcal U}}

\newcommand{\sK}{{\mathbb K}}
\newcommand{\sL}{{\mathbb L}}

\newcommand{\sT}{{\mathbb T}}

\newcommand{\tB}{{\widetilde{B}}}

\newcommand{\tI}{{\widetilde{I}}}

\newcommand{\tL}{{\widetilde{L}}}
\newcommand{\tM}{{\widetilde{M}}}

\newcommand{\tP}{{\widetilde{P}}}
\newcommand{\tQ}{{\widetilde{Q}}}

\newcommand{\tU}{{\widetilde{U}}}
\newcommand{\tV}{{\widetilde{V}}}
\newcommand{\tW}{{\widetilde{W}}}

\newcommand{\tb}{{\widetilde{b}}}

\newcommand{\tg}{{\tilde{g}}}

\begin{document}
\title{Quantum Cluster Variables via Serre Polynomials}
\author{Fan QIN}
\address[]{Fan Qin, Universit\'{e} Paris Diderot - Paris 7, Institut de Math\'{e}\-ma\-ti\-que de Jussieu, UMR 7586 du CNRS, 175 rue du chevaleret, 75013, Paris, France}
\email{qinfan@math.jussieu.fr}%add

\begin{abstract}
%% Text of abstract
For skew-symmetric acyclic quantum cluster algebras, we express the quantum $F$-polynomials and the quantum cluster monomials in terms of Serre polynomials of quiver Grassmannians of rigid modules. As byproducts, we obtain the existence of counting polynomials for these varieties and the positivity conjecture with respect to acyclic seeds. These results complete previous work by Caldero and Reineke and confirm a recent conjecture by Rupel.
\end{abstract}

\thanks{The author is financially supported by a CSC scholarship.}

\maketitle

\tableofcontents
%%%%%%%%%%%%%%%%%%%%%%%%%%
%%			Main content
%%%%%%%%%%%%%%%%%%%%%%%%%%

\section{Introduction}\label{sec:intro}

\subsection{Motivation}

Cluster algebras were invented by Fomin and Ze\-le\-vin\-sky \cite{FominZelevinsky02} in order to provide a combinatorial approach to canonical bases and total positivity. They are commutative subalgebras of Laurent polynomial rings. Later, quantum cluster algebras were introduced in \cite{BerensteinZelevinsky05}. Therefore, since its very beginning, the theory of (commutative and quantum) cluster algebras is closely related to topics in quantum groups and canonical bases.

Let us point out another important aspect of quantum cluster theory. Recall that in this theory a recursive operation called mutation plays the central role. In \cite{KontsevichSoibelman08} Kontsevich and Soibelman interpreted mutations of quantum cluster algebras in their study of Donaldson-Thomas invariants of certain triangulated categories of Calabi-Yau dimension~$3$. This establishes a remarkable link between quantum cluster theory and triangulated categories.

Moreover, the occurrence of quantum cluster algebras precedes that of commutative cluster algebras in the study of Donaldson-Thomas invariants. Therefore, quantum cluster algebras are important and interesting despite their combinatorial similarity with their commutative counterparts.

\subsection{Contents}
In the present paper, we keep the role of triangulated categories to a minimum, and restrict our attention to quantum cluster variables, which are certain generators of quantum cluster algebras defined recursively by mutations. We refer to Section \ref{sec:main} for the necessary definitions and a more detailed summary of the main results.

\subsubsection*{The main theorem: A refined CC-formula via Serre polynomials}

Recall that the commutative cluster variables can be expressed in terms of the Euler characteristics of quiver Grassmannians \cite{CalderoChapoton06}. The expression is given by what we call the (numerical) CC-formula. Then naturally we ask for a refined CC-formula to compute quantum cluster variables. Notice that the coefficients of quantum cluster variables are Laurent polynomials in $q^{\Hf }$, \st
\begin{enumerate}
\item they specialize to the Euler characteristics of quiver Grassmannians under the quasi-classical limit $q^\Hf\ra 1$, and
\item they are invariant under the involution $q^\Hf\mapsto q^{-\Hf}$.
\end{enumerate}
In addition, for an acyclic quiver, the quiver Grassmannians associated to commutative cluster variables are smooth \cite{CalderoReineke08}. All these facts remind us of the cohomology of compact smooth varieties. Accordingly, under the assumption that the quiver is acyclic, we propose a refined CC-formula via Serre polynomials as below. 

Let $\tQ$ be an ice quiver, \ie a quiver with frozen vertices. Assume that the matrix of $\tQ$ can be completed into a unitally compatible pair, to which we associate a quantum cluster algebra $\cA^q$. Further suppose that the non frozen part $Q$ of $\tQ$ is acyclic. Endow $\tQ$ with a generic potential in the sense of \cite{DerksenWeymanZelevinsky08} and consider the associated presentable cluster category $\cD$ as introduced in \cite{Plamondon10a}.
\begin{Def*}[Definition \ref{def:X_M}, refined CC-formula]
For any coefficient-free rigid object $M$ of $\cD$, we denote by $m$ the class of $\Ext^1_{\cC}(T,M)$ in $\Kz(\mod kQ)$, and associate to $M$ the following element in the quantum torus $\cT$:
\begin{align*}
X_M=\sum_e E(\Gr_e (\Ext^1_{\cC}(T,M))) q^{-\Hf \langle e,m-e \rangle} X^{\ind_T(M)-\phi(e)}.
\end{align*}
\end{Def*}
Here the symbol $E(\Gr_e (\Ext^1_{\cC}(T,M)))$ denotes the Serre polynomial of $\Gr_e (\Ext^1_{\cC}(T,M))$. This variety is smooth and projective (because $Q$ is acyclic and $M$ is rigid) so that we have
\[
E(\Gr_e (\Ext^1_{\cC}(T,M)))=\sum_i (-1)^i \dim H^i(\Gr_e (\Ext^1_{\cC}(T,M)))q^{\Hf i}.
\]
The main theorem (Theorem \ref{thm:X_M}) claims that all the quantum cluster variables, and further all the quantum cluster monomials of $\cA^q$, take this form.

\subsubsection*{The ingredients of the proof}

Inspired by Katz's Theorem \ref{thm:Katz_thm}, which links Serre polynomials to counting polynomials, we follow an indirect approach to prove the main theorem, which roughly consists of the following steps. 
\begin{itemize}
\item In the definition of the quantum cluster algebra $\cA^q$, we replace the formal parameter $q^\Hf$ by the number $|\sK|^\Hf$ for each finite field $\sK$. Then we obtain what we call the specialized quantum cluster algebra $\cA^\sK$. We further notice that there is a map $ev_\sK$ which sends $q^\Hf$ to $|\sK|^\Hf$ and the quantum cluster variables of $\cA^q$ to the specialized quantum cluster variables of $\cA^\sK$.
\item
We prove an analogue of the main theorem for $\cA^\sK$, \ie we show that the specialized quantum cluster variables of $\cA^\sK$ can be expressed in terms of the numbers of $\sK$-points of quiver Grassmannians. To see this, we observe that these numbers equal sums of Hall numbers, and then use Green's Theorem \cite{Green95} to complete the proof following an idea of Hubery \cite{Hubery}. As a consequence, these quiver Grassmannians have counting polynomials, which agree with the coefficients of quantum cluster variables up to $q^\Hf$-power factors.
\item
Theorem \ref{thm:Katz_thm} claims that if a scheme has a counting polynomial, then it must be the Serre polynomial. So the coefficients of quantum cluster variables agree with the Serre polynomials up to $q^\Hf$-power factors.
\end{itemize}
\subsubsection*{By-products}
As by-products of the main theorem, for an acyclic quiver $Q$, we obtain that 
\begin{enumerate}
\item
the quantum $F$-polynomials can be expressed in terms of Serre polynomials, and that
\item
the submodule Grassmannians of rigid modules have their cohomology concentrated in even degrees, and have their Serre polynomials as counting polynomials. 
\end{enumerate}

The above geometric properties of submodule Grassmannians first appeared in a theorem of \cite{CalderoReineke08}, though the proof there contains a gap as pointed out by \cite{Nakajima09}. It is somehow surprising that in the present paper we obtain a proof via quantum cluster algebras. This might be viewed as a consequence of the link between quantum cluster algebras and algebraic geometry. 

While the author was finishing the present paper, Rupel \cite{Rupel10} posed a conjecture that quantum cluster variables could be expressed in terms of counting polynomials of quiver Grassmannians, and proved it for those variables which lie in almost acyclic clusters using combinatorial methods. Notice that this conjecture already appeared in embryonic form in the introduction of \cite{CalderoChapoton06}. Clearly, at least in the case of non valued quivers, this conjecture is confirmed by our main theorem, since the Serre polynomials are also the counting polynomials. Demonet's work \cite{Demonet10} suggests that our techniques may also be useful in the valued case.

After the first version of this article was posted on Arxiv, Professor Hiraku Nakajima kindly informed the author about a geometric proof \cite{NakajimaEmail} for the vanishing of the odd cohomology of quiver Grassmannians. That proof is based on the Fourier-Sato-Deligne transform and the celebrated Decomposition theorem. Details will appear in an updated version of \cite{Nakajima09}.

\section*{Acknowledgments}\label{sec:ack}

The author would like to express his sincere thanks to his thesis advisor Bernhard Keller for all the interesting discussions and a lot of other help. He also thanks Dong Yang for many useful discussions in representation theory. He thanks Pierre-Guy Plamondon for sharing and discussing his preprints \cite{Plamondon10a} \cite{Plamondon10b}. And he thanks Tamas Hausel for his inspiring series of lectures at CIRM in September 2009, where he learned about the wonderful properties of Serre polynomials. The author thanks Professor Hiraku Nakajima for informing him about his geometric proof of the vanishing of odd cohomology of quiver Grassmannians.

\section{Preliminaries}\label{sec:reminders}
	\subsection{Quantum cluster algebras}\label{sec:quantum_cluster_algebra}
		\subsubsection*{Generalized quantum cluster algebras}
		
Let $R$ be an integral domain (\ie an integral commutative ring) and let $v$ be an invertible element in $R$. We shall define generalized quantum cluster algebras over $(R,v)$, in analogy with quantum cluster algebras as defined in \cite{BerensteinZelevinsky05}. In fact, in the present paper, we shall only be interested in the following cases:
\begin{itemize}
	\item When $(R,v)=(\Z, 1)$, we obtain the usual commutative cluster algebras $\cA^\Z$.
	\item When $(R,v)=(\Z[q^{\pm\Hf }],q^{\Hf })$, where $q^\Hf$ is a formal parameter, we obtain quantum cluster algebras $\cA^q$.
	\item When $(R,v)=(\Z[|\sK|^{\pm\Hf }],|\sK|^{\Hf })$, where $\sK$ is a finite field, we obtain what we call specialized quantum cluster algebras $\cA^\sK$.
\end{itemize}

Let $m\geq n$ be two positive integers. Let $\Lambda$ be an $m\times m$ skew-symmetric integer matrix and $\tB$ an $m\times n$ integer matrix. The upper $n\times n$ submatrix of $\tB$, denoted by $B$, is called the \emph{principal part} of $\tB$.
\begin{Def}[Compatible pair]\label{def:compatible}
The pair $(\Lambda, \tB)$ is called \emph{compatible }if we have
\begin{align}\label{eq:BZ_compatible}
\Lambda(-\tB)=\begin{bmatrix}D\\0 \end{bmatrix}
\end{align}
for some $n\times n$ diagonal matrix $D$ whose diagonal entries are strictly positive integers. It is called a \emph{unitally} compatible pair if moreover $D$ is the identity matrix.
\end{Def}

Let $(\Lambda,\tB)$ be a compatible pair. The component $\Lambda$ is called the \emph{$\Lambda$-matrix} of $(\Lambda,\tB)$, and the component $\tB$ the \emph{$B$-matrix} of $(\Lambda,\tB)$.

\begin{Prop}\cite[Proposition 3.3]{BerensteinZelevinsky05}
The $B$-matrix $\tB$ has full rank $n$, and the product $D B$ is skew-symmetric.
\end{Prop}

Let $\aL$ be a rank $m$ lattice over $\Z$ and $\lambda$ a skew-symmetric bilinear form on $L$. 
\begin{Def}[Quantum torus]\label{def:quantum_torus}
The \emph{quantum torus }$\cT=\cT(L,\lambda)$ over $(R,v)$ is the $R$-algebra generated by the symbols $X^g$, $g\in L$, subject to the relations
\[
X^{g}X^{h}=v^{\lambda(g,h)}X^{g+h}
\]
for $g$ and $h$ in $L$.
\end{Def}
Denote the map from $L$ to $\cT$ sending $g$ to $X^g$ by $X_{\cT}$. Notice that the torus $\cT$ is contained in its \emph{skew-field of fractions}, which is denoted by $\cF$ \cite[Appendix]{BerensteinZelevinsky05}. In the following statements, we identify the lattice $\Z^m$ with $L$ by choosing a basis for $L$.

\begin{Def}[Toric frame]
A \emph{toric frame }is a map $X:\Z^m\longrightarrow \cF$, such that for some automorphism $\rho$ of the skew-field $\cF$, and some automorphism $\eta$ of the lattice $\Z^m$, one has $X(c)=\rho(X^{\eta(c)})$ for all $c$ in $\Z^m$.
\end{Def}

Denote the operation of matrix transposition by $(\ )^T$.
\begin{Def}[Seed]
A \emph{seed} is a triple $(\Lambda,\tB, X)$, \st
\begin{enumerate}
\item $(\Lambda,\tB)$ is a compatible pair,
\item $X$ is a toric frame, and
\item $X(g)X(h)=v^{g^T \Lambda h}X(g+h)$ for all $g$, $h$ in $\Z^m$.
\end{enumerate}
It is said to be unitally if $(\Lambda, \tB)$ is.
\end{Def}
We often write $\Lambda(g,h)$ for $g^T \Lambda h$.

Let $(\Lambda,\tB, X)$ be a seed. We associate to it the following elements in $\cT$.
\begin{Def}[Quantum Cluster variables]\label{def:cluster_variable}
The \emph{$X$-variables} of the seed $(\Lambda,\tB, X)$ are defined as $X_i=X(e_i), 1\leq i \leq m$. The set $\{X_k, 1\leq k\leq n\}$ is called the \emph{cluster}; its elements are called \emph{quantum cluster variables}. The \emph{quantum cluster monomials }are those $X(e)$ \st $e_i\geq 0$ for $1\leq i\leq n$ and $e_j=0$ for $j>n$. 
\end{Def}

A \emph{sign} $\epsilon$ is an element in $\{-1,+1\}$. Let $b$ be a vector in $\Z^m$ which does not belong to the kernel of the form $\Lambda(\ ,\ )$. Then the image of the map taking $g\in\Z^m$ to $\Lambda(b,g)$ is a nonzero ideal of $\Z$. Let $d(b)$ denote its positive generator.
\begin{Prop}[Automorphism]\cite[Proposition 4.2]{BerensteinZelevinsky05}
For any $1\leq k\leq n$ and any sign $\epsilon$, the skew-field $\cF$ has an automorphism $\rho_{b, \epsilon}$ \st
\begin{align*}
\rho_{b, \epsilon}(X(c))=(P^{|d|}_{-\epsilon_d b,-\epsilon_d \epsilon})^{-\epsilon_d} X(c)
\end{align*}
for $c\in\Z^m$, where $d=\frac{\Lambda(b,c)}{d(b)}$, $\epsilon_d=\sign d$ is the sign of $d$, and for $r\in\N$ we put
\[
P^{r}_{b, \epsilon}=\prod^{r}_{p=1} (1+v^{\epsilon (2p-1)d(b)}X(b\epsilon)).
\]
\end{Prop}

Denote by $b_{ij}$ the entry in position $(i,j)$ of $\tB$. For any $1\leq k\leq n$ and any sign $\epsilon$, we associate to $\tB$ an $m\times m$ matrix $E_\epsilon$ whose entry in position $(i,j)$ is 
\begin{align}\label{eq:E_epsilon}
e_{ij}=\left\{
\begin{array}{ll}
\delta_{ij} & \textrm{if $j\neq k$}\\
-1 & \textrm{if $i=j=k$}\\
max(0,-\epsilon b_{ik}) & \textrm{if $i\neq k, j=k$},
\end{array} \right.
\end{align}
and an $n\times n$ matrix $F_\epsilon$ whose entry in position $(i,j)$ is
\[
f_{ij}=\left\{
\begin{array}{ll}
\delta_{ij} & \textrm{if $i\neq k$}\\
-1 & \textrm{if $i=j=k$}\\
max(0,\epsilon b_{kj} ) & \textrm{if $i=k$, $j\neq k$}.
\end{array} \right.
\]

\begin{Def}[Seed Mutation]\label{def:seed_mutation}
For any $1\leq k\leq n$ and any sign $\epsilon$, we define the mutated seed $(\Lambda',\tB', X')=\mu_k(\Lambda,\tB, X)$ by
\[
(\Lambda',\tB')=(E_\epsilon ^T \Lambda  E_\epsilon , E_\epsilon  \tB  F_\epsilon ),
\]
and 
\begin{align}\label{eq:frame_mutation}
X'(c)=\rho_{\tb^k , \epsilon}(X(E_\epsilon  c))
\end{align}
for all $c$ in $\Z^m$. Here $\tb^k $ is the $k$-th column of $\tB  $.
\end{Def}
One should notice that $\mu_k$ is an involution, and is independent of the choice of $\epsilon$.

Obviously, the map $X'$ in the above definition is determined by its values on the natural basis $\{e_i\}$ of $\Z^m$. For an integer $x$, we write $[x]_+$ for $max(0,x)$.
\begin{Lem}[Mutation of quantum cluster variables]\label{lem:cluster_mutation}
The toric frame $X'$ in Definition \ref{def:seed_mutation} is determined by
\begin{align}
X_k'&=X (\sum_{1\leq i\leq m}[b_{ik}]_{+} e_i -e_k)+X (\sum_{1\leq j\leq m}[-b_{jk}]_{+} e_j -e_k),\\
X_k'&=X_i ,\quad 1\leq i\leq m,\quad i\neq k.
\end{align}
\end{Lem}

\begin{Cor} \label{cor:cluster_mutation}
We have
\begin{align}
\begin{split}
X_k X_k'=&v^{\Lambda(e_k,\sum_{1\leq i\leq m}[b_{ik}]_{+} e_i)}X (\sum_{1\leq i\leq m}[b_{ik}]_{+} e_i)\\
&\qquad+v^{\Lambda(e_k,\sum_{1\leq j\leq m}[-b_{jk}]_{+} e_j)}X (\sum_{1\leq j\leq m}[-b_{jk}]_{+} e_j).
\end{split}
\end{align}
\end{Cor}

\begin{Lem}[Involution]
Let $X$ be a toric frame. Assume that $R$ is endowed with an involution sending each element $r$ to $\bar{r}$, \st $\bar{v}$ equals $v^{-1}$.
\begin{enumerate}
\item There is a unique involutive antiautomorphism $\sigma_X$ of the ring $\cT$ \st \linebreak[4]$\sigma_X(X(c))=X(c)$ for all $c$ in $\Z^m$ and $\sigma_X(r)=\overline{r}$ for all $r$ in $R$.
\item If $(\Lambda' ,\tB'  , X' )$ is a seed obtained from a seed $(\Lambda ,\tB  , X )$ by a mutation, then the map $\sigma_{X'}$ equals $\sigma_X$. In particular, we have $\sigma_X(X'(c))=X'(c)$ for all $c$ in $\Z^m$.
\end{enumerate}
\end{Lem}

Let $\cT=\cT(\Z^m, \Lambda)$ be a quantum torus over $(R,v)$ with skew-field of fractions $\cF$ and $(\Lambda, \tB)$ a compatible pair. Let $\sT_n$ be an $n$-regular tree with root $t_0$. There is a unique way of associating a seed $(\Lambda(t),\tB(t),X(t))$ with each vertex $t$ of $\sT_n$ \st we have
\begin{enumerate}
\item $(\Lambda(t_0),\tB(t_0),X(t_0))=(\Lambda,\tB,X_\cT)$, and
\item if two vertices $t$ and $t'$ are linked by an edge labeled $k$, then the seed \linebreak[4]$(\Lambda(t'),\tB(t'),X(t'))$ is obtained from $(\Lambda(t),\tB(t),X(t))$ by the mutation at $k$.
\end{enumerate}
Notice that for $j>n$, the elements $X_{j}(t)$ do not depend on $t$.
\begin{Def}[Quantum cluster algebra over $(R,v)$]
The \emph{quantum cluster algebra} $\cA$ over $(R,v)$ is the $R$-subalgebra of $\cF$ generated by the quantum cluster variables $X_i(t)$ for all the vertices $t$ of $\sT_n$ and $1\leq i\leq n$, and the elements $X_{j}(t_0)$ and $(X_{j}(t_0))^{-1}$ for all $j>n$.
\end{Def}

\begin{Thm}[Laurent phenomenon]\cite[Section 5]{BerensteinZelevinsky05}\label{thm:Laurent}
The quantum cluster algebra $\cA$ is a subalgebra of $\cT$.
\end{Thm}
		
\subsubsection*{Quantum \texorpdfstring{$F$}{F}-polynomials and extended g-vectors}
		
		We recall and adapt some results in \cite{Tran09} to our setting.

\begin{Def}
\cite[Definition 3.7]{Tran09} Let $\cA^\Z$ be the commutative cluster algebra associated with $\tB$ as in \cite{FominZelevinsky07}. For each vertex $t$ of $\sT_n$ and each $1\leq l\leq n$, the \emph{extended $g$-vector} associated to $t$ and $l$ is defined as the vector $\tg_l(t)\in\Z^m$, \st the cluster variable associated to $t$ and $l$ is
\begin{align}
X_l^\Z(t)=F_l^\Z(t)|_{y^e\mapsto X^{\tB e}}X_l^{\tg_l(t)},
\end{align}
where $F_l^\Z(t)$ is the commutative $F$-polynomial associated to $t$ and $l$ as in \cite{FominZelevinsky07}.
\end{Def}

%According to \cite{FominZelevinsky07}, we have $g$-vector $\tg_l(t)=(g_l(t),g_l^c(t))$, where %$g_l(t)\in\Z^n$ is the usual $g$-vector, and complementary part $g_l^c(t)$ can be computed by specializing %$F_l(t)$ over some tropcial semi-field.
\begin{Rem}
Assume that the principal part $B$ of $\tB$ is skew-symmetric. We can find an $m\times m$ skew-symmetric matrix $\tB^\circ$ whose left $n$ columns form the matrix $\tB$. Then the commutative cluster algebra $\cA^\Z_{\tB^\circ}$ associated to $\tB^\circ$ contains $\cA^\Z$ as a subalgebra. As a consequence, the extended $g$-vectors $\tg_l(t)$ of $\cA^\Z$ equal the usual $g$-vectors $g_l(t)$ of $\cA^\Z_{\tB^\circ}$.
\end{Rem}

Let $B$ be an $n\times n$ integer matrix and $D$ an $n\times n$ diagonal integer matrix \st the product $DB$ is antisymmetric. 
\begin{Def}\label{def:torus_matrix}
The \emph{quantum torus} $\cR(DB)$ associated to the antisymmetric matrix $DB$ is the $\Rq$-algebra generated by the symbols $y^e$, $e\in \Z^n$, subject to the relations
\[
y^e y^f=q^{\Hf e^T DB f}y^{e+f}
\]
for $e,f\in\Z^n$.
\end{Def}

\begin{Thm}\label{thm:Tran_thm}
\cite[Theorem 5.3, Theorem 6.1, Remark 6.3]{Tran09} Let $\cA^q$ be the quantum cluster algebra associated with a compatible pair $(\Lambda, \tB)$ \st the principal part $B$ of $\tB$ is antisymmetric. Let $D$ be as in \eqref{eq:BZ_compatible}. For any vertex $t$ of $\sT_n$ and any $1\leq l\leq n$, there exists an element $F_l(t)$ of $\cR(DB)$ such that we have
\begin{align}
X_l(t)=F_l(t)|_{y^e\mapsto X^{\tB e}}X^{\tg_l(t)}.
\end{align}
The element $F_l(t)$ is called the \emph{quantum $F$-polynomial} associated to $t$ and $l$. It specializes to $F_l^\Z(t)$ under the quasi-classical limit $q^{\Hf}\ra 1$.

Moreover, the quantum $F$-polynomial $F_l(t)$ only depends on $l$, $t$, $B$ and $D$ (and not on the non principal part of $\tB$).
\end{Thm}
		
\subsection{Specialization maps}\label{sec:specialization}
		
		Let $L$ be a rank $m$ $\Z$-lattice endowed with a skew-symmetric bilinear form $\lambda$. Associate to $(L,\lambda)$ the quantum torus $\cT$ over $(\Rq,q^\Hf)$, the quantum torus $\cT^\Z$ over $(\Z,1)$, and for each finite field $\sK$ the quantum torus $\cT^\sK$ over $(\Rk,|\sK|^\Hf)$. Notice that the torus $\cT^\Z$ is just the Laurent polynomial ring in $m$ variables.
		
		Clearly, there are unique algebra homomorphisms $ev_1:\cT\ra\cT^\Z$ and $ev_\sK:\cT\ra\cT^\sK$ which send $X^g$ to $X^g$ for $g$ in $L$ and $v$ to $1$ respectively $|\sK|^\Hf$. It is easy to see that $ev_1$ induces an isomorphism $\cT/(q^\Hf-1)\iso \cT^\Z$ and $ev_\sK$ an isomorphism $\cT/(q^\Hf-|\sK|^\Hf)\iso \cT^\sK$. Therefore we call $ev_1$ and $ev_\sK$ \emph{specialization maps}. 
		
		To a compatible pair $(\Lambda, \tB)$, associate the \emph{quantum cluster algebra} $\cA^q$ over $(\Rq,q^\Hf)$ with the initial seed $(\Lambda, \tB, X_{\cT})$, the \emph{commutative cluster algebra} $\cA^\Z$ over $(\Z,1)$ with the initial seed $(\Lambda, \tB, X_{\cT^\Z})$, and the \emph{specialized quantum cluster algebra} $\cA^\sK$ over $(\Rk$,\linebreak[4]$|\sK|^\Hf)$ with the initial seed $(\Lambda, \tB, X_{\cT^\sK})$. For any $1\leq i\leq m$ and any vertex $t$, denote the associated (quantum) $X$-variables of $\cA^q$, $\cA^\Z$ and $\cA^\sK$ by $X_i(t)$, $X_i^\Z(t)$, $X_i^\sK(t)$ respectively.

We have seen in Theorem \ref{thm:Laurent} that $\cA^q\subset \cT$, $\cA^\Z\subset \cT^\Z$, $\cA^\sK\subset \cT^\sK$. Since quantum tori do not have zero divisors, inductively applying Corollary \ref{cor:cluster_mutation}, we obtain that $ev_1(X_i(t))=X_i^\Z(t)$ and $ev_\sK(X_i(t))=X_i^\sK(t)$ for any $1\leq i\leq m$ and any vertex $t$.

Recall that the $B$-matrix $\tB$ is of full rank. Define a partial order on $\Z^m$ as follows.
\begin{Def}\label{def:partial_order}
For any elements $g$ and $h$ in $\Z^m$, $g$ is said to be \emph{strictly larger }than $h$ if there exists some nonzero element $e$ in $\Z^n$ with non-negative coordinates \st we have $g=h+\tB e$.
\end{Def}

Let $U=\sum_{g\in \Z^m} P_g X^g$, $P_g\in R$, be an element in a quantum torus over $(R,v)$ as in Definition \ref{def:quantum_torus}. The \emph{support} of $U$ is the set of elements $g$ of $\Z^m$ \st $P_g$ is nonzero. A \emph{minimal degree} of $U$ is an element of its support which is minimal with respect to the partial order in the above definition.

\begin{Lem}\label{lem:multi_degree}
For any elements $U$ and $U'$ of a quantum torus over $(R,v)$, if they both have a unique minimal degree which we denote by $\deg U$ and $\deg U'$ respectively, then the product $U\cdot U'$ has the unique minimal degree $\deg U+\deg U'$.
\end{Lem}

\begin{Prop}\label{prop:spec_bijective}
(1) The specialization map $ev_1$ induces a support-preserving bijection from the set of (quantum) cluster monomials of $\cA^q$ to that of $\cA^\Z$.

(2) The specialization map $ev_\sK$ induces a support-preserving bijection from the set of quantum cluster monomials of $\cA^q$ to that of $\cA^\sK$.
\end{Prop}
\begin{proof}
(1) We know that the map $ev_1$ is surjective. It remains to show its injectivity.

For this, we adopt the idea of the proof of Theorem 6.1 of \cite{BerensteinZelevinsky05}. Fix two vertices $t$ and $t'$ of $\sT_n$ and two quantum cluster monomials  $X(t)(e)$ and $X(t')(e')$ as in Definition \ref{def:cluster_variable} \st they have the same image under $ev_1$. Consider the seeds $(\Lambda(t),\tB(t),X(t))$ and $(\Lambda(t'),\tB(t'),X(t'))$ associated to $t$ and $t'$ respectively. By Section 6 of \cite{BerensteinZelevinsky05}, every quantum cluster variable $X_j(t')$, $1\leq j\leq n$, has an expansion of the form $X_j(t')=\sum_{d\in\Z^m} J_d X(t)(d)$, where the coefficients $J_d$ belong to $\Rq$ and can also be written as subtraction-free rational expressions in $q^\Hf$. Then this property also holds for $X(t')(e')$, and we can write 
\[
X(t')(e')=\sum_{d\in\Z^m} P_d X(t)(d)
\]
where the coefficients $P_d$ can be written as subtraction-free rational expressions. Then for any nonzero coefficient $P_d$, we have $ev_1(P_d)\neq 0$. Now since $ev_1(X(t')(e'))=ev_1(X(t)(e))$ and since, for fixed $t$, the monomials $X^\Z(t)(d)$, $d\in \Z^m$, are linearly independent over $\Q$, we have $X(t')(e')=p X(t)(e)$ for some $p$ in $\Rq$. Because of the symmetry between $X(t')(e')$ and $X(t)(e)$, $p$ is invertible. Finally, since quantum cluster monomials are involution invariant, we obtain that $p$ equals $1$.

(2) The proof is the similar to that of (1).
\end{proof}

Assume that the principal part $B$ of $\tB$ is skew-symmetric. We have the following result.
\begin{Prop}\label{prop:minimal_degree}
Every quantum cluster monomial of $\cA^\sK$ has a unique minimal degree. Moreover, different quantum cluster monomials of $\cA^\sK$ have different minimal degrees.
\end{Prop}
\begin{proof}
Each commutative $F$-polynomial $F^\Z_i(t)$ has a non-zero constant term, \cf\cite{DerksenWeymanZelevinsky09}. Therefore the quantum $F$-polynomial $F_i(t)$ also has a non-zero constant term. By Theorem \ref{thm:CC_formula} below and Theorem \ref{thm:Tran_thm}, we see that the cluster variable $X^\Z_i(t)$ and the quantum cluster variable $X_i(t)$ have the same unique minimal degree $\ind_T T_i(t)$.

Inductively applying Lemma \ref{lem:multi_degree} which claims that minimal degrees are additive, we obtain that every quantum cluster monomial $X(t)(e)$ as in Definition \ref{def:cluster_variable} and its image $ev_\sK(X(t)(e))=X^\sK(t)(e)$ have the same minimal degree $\sum_{1\leq i\leq n}e_i \ind_T T_i(t)=\ind_T T(t)(e)$, where we write $T(t)(e)$ for the object $\oplus_i T_i(t)^{\oplus e_i}$. By \cite[Theorem 3.7]{Plamondon10b}, different coefficient-free rigid objects have different indices. In addition, by Theorem \ref{thm:CC_formula} below and Proposition \ref{prop:spec_bijective} we have a bijection between the set of coefficient-free rigid objects and the set of quantum cluster monomials of $\cA^\sK$, which sends an object $T(t)(e)$ to the quantum cluster monomial $X^\sK(t)(e)$. Therefore, different quantum cluster monomials of $\cA^\sK$ have different minimal degrees.
\end{proof}

%This lemma implies the following result.
%\begin{Prop}\label{prop:specialization}
%The specialization maps $ev_1$ and $ev_\sK$ induce algebra homorphisms between quantum cluster algebras.

%Moreover, for any vertex $t$ of $\sT_n$ and any $1\leq i\leq m$, the map $ev_1$ sends $X_i(t)$ to $X_i^\Z(t)$ and induces a bijection from the set of quantum cluster monomials of $\cA^q$ to the set of cluster monomials of $\cA^\Z$, and the map $ev_\sK$ sends $X_i(t)$ to $X_i^\sK(t)$ and induces a bijection from the set of quantum cluster monomials of $\cA^q$ to the set of quantum cluster monomials of $\cA^\sK$.
%\end{Prop}
%We denote the inverse maps of the bijections by $ev_1^{-1}$ and $ev_\sK^{-1}$ respectively.

\subsection{The CC-formula for commutative cluster variables} \label{sec:cluster_category}

	Let $m\geq n$ be two positive integers and $\tQ$ a quiver without loops or $2$-cycles and with vertex set $\{1,\ldots,m\}$. Denote the subset $\{n+1,\dots,m\}$ by $C$. The elements in $C$ are called the \emph{frozen vertices }, and $\tQ$ is called an \emph{ice quiver}. The full subquiver $Q$ on the vertices $1,\ldots,n$ is called the \emph{principal part} of $\tQ$. Let $Q_0$ denote the set $\{1,\ldots,n\}$ and $\tQ_0$ denote the set $\{1,\ldots,m\}$.
	
Let $\tB$ be the $m\times n$ matrix associated to the ice quiver $\tQ$, \ie its entry in position $(i,j)$ is
\[
b_{ij}=\sharp\{\mathrm{arrows}\, i\ra j\}-\sharp\{\mathrm{arrows}\, j\ra i\}
\]
for $1\leq i\leq m$, $1\leq j\leq n$. We associate to $\tB$ the cluster algebra $\cA^\Z$ as in \cite{FominZelevinsky07}. 
	
Let the base field $k$ be the complex field $\C$.	Let $\tW$ be a generic potential on $\tQ$ in the sense of \cite{DerksenWeymanZelevinsky08}. As in \cite{KellerYang09}, with the quiver with potential $(\tQ,\tW)$ we can associate the Ginzburg algebra $\Gamma=\Gamma(\tQ,\tW)$. Denote the perfect derived category of $\Gamma$ by $\per\Gamma$ and denote the full subcategory of $\per\Gamma$ whose objects are dg modules with finite dimensional homology by $\Dfd\Gamma$. The \emph{generalized cluster category} $\cC=\cC_{(\tQ,\tW)}$ in the sense of \cite{Amiot09} is the quotient category
\[
\cC=\per\Gamma/\Dfd\Gamma.
\]
Denote the quotient functor by $\pi:\per\Gamma\ra \cC$ and define
\begin{align*}
T_i&=\pi(e_i\Gamma),\quad 1\leq i\leq m\\
T&=\oplus_{1\leq i\leq m} T_i.
\end{align*}
It is shown in \cite{Plamondon10c} that the endomorphism algebra of $T$ is isomorphic to $H^0\Gamma$.

For any triangulated category $\cU$ and any rigid object $X$ of $\cU$, we define the subcategory $\pr_\cU (X)$ of $\cU$ to be the full subcategory consisting of the objects $M$ \st there exists a triangle in $\cU$
\[
M_1\ra M_0\ra M\ra \Sigma M_1,
\]
for some $M_1$ and $M_0$ in $\add X$. The \emph{presentable cluster category} $\cD\subset \cC$ is defined as the full subcategory consisting of the objects $M$ \st
\[
M\in\pr_\cC (T)\cap\pr_\cC (\Sigma^{-1}T)\quad\text{and}\quad\dim\Ext^1_\cC(T,M)<\infty,
\]
\cf \cite{Plamondon10a}.

We refer to \cite{Plamondon10a} for the definition of the iterated mutations of the object $T$. There is a unique way of associating an object $T(t)=\oplus_{1\leq i\leq m} T_i(t)$ of $\cD$ with each vertex $t$ of $\sT_n$ \st we have
\begin{enumerate}
\item $T(t_0)=T$, and
\item if two vertices $t$ and $t'$ are linked by an edge labeled $k$, then the object $T(t')$ is obtained from $T(t)$ by the mutation at $k$. 
\end{enumerate}

Let $\cF\subset\per\Gamma$ denote the full subcategory $\pr_{\per\Gamma}(\Gamma)$. The quotient functor $\pi:\per\Gamma\ra \cC$ induces an equivalence $\cF\iso\pr_\cC(T)$. Denote by $\pi^{-1}$ the inverse equivalence. For an object $M\in \pr_{\cC}(T)$, we define its \emph{index }$\ind_T M$ as the class $[\pi^{-1}M]$ in $\Kz(\per\Gamma)$. 

\begin{Thm}\cite{Plamondon10a}\label{thm:g_vector_coordinate}
(1) For any vertex $t$ of $\sT_n$, the classes $[\ind_T T_i(t)]$ form a basis of $\Kz(\per\Gamma)$. 

(2) For a class $[P]$ in $\Kz(\per\Gamma)$, let $[[P]:T_i(t)]$ denote its $i$th coordinate in this basis. Then we have $([\ind_T T_i(t):T_j])_{1\leq j\leq m}=\tg_i(t)$.
\end{Thm}

\begin{Def}[Coefficient-free objects]\label{def:CoefFree}
An object $M$ in $\cC$ is called \emph{coefficient-free} if
\begin{enumerate}
\item  the object $M$ does not contain a direct summand $T_i$, $i>n$, and
\item  the space $\Ext^1_{\cC}(T_i,M)$ vanishes for $i>n$.
\end{enumerate}
\end{Def}
For a coefficient-free object $M\in\cC$, the space $\Ext^1_{\cC}(T,M)$ is a right $H^0\Gamma$-module whose support is concentrated on $Q$. Thus, it can be viewed as a $\cP(Q,W)$-module, where $W$ is the potential on $Q$ obtained from $\tW$ by deleting all cycles through vertices $j>n$ and $\cP(Q,W)$ is the Jacobi algebra of $(Q,W)$. Denote by $\phi:\Kz(\mod \cP(Q,W))\ra \Kz(\per\Gamma)$ the map induced by the composition of inclusions $\mod\cP(Q,W)\ra \Dfd\Gamma \ra \per\Gamma$. For any vertex $i$ of $\tQ$, it sends $[S_i]$ to
\[
\sum_{\mathrm{arrows}\, i\ra j}[e_j\Gamma]-\sum_{\mathrm{arrows}\, l\ra i}[e_l\Gamma],
\]
as one easily checks using the minimal cofibrant resolution of the simple dg $\Gamma$-module $S_i$, \cf \cite{KellerYang09}.

\begin{Def}[CC-formula]\label{def:CC_formula}
For any coefficient-free and rigid object $M\in\cD$, we denote by $m$ the class of $\Ext^1_{\cC}(T,M)$ in $\Kz(\mod kQ)$, and associate to $M$ the following element in $\cT^\Z$:
\begin{align}\label{eq:CC_formula}
X^\Z_M=\sum_e \chi(\Gr_e (\Ext^1_{\cC}(T,M))) X^{\ind_T(M)-\phi( e)}.
\end{align}
where $\Gr_e (\Ext^1_{\cC}(T,M))$ is the submodule Grassmannian of $\Ext^1_{\cC}(T,M)$ whose $\C$-points are the submodules of the class $e$ in $\Kz(\mod\cP(Q,W))$.
\end{Def}

The following theorem is a consequence of the results in \cite{Plamondon10c} and \cite{CalderoKeller06}.
\begin{Thm} \label{thm:CC_formula}
Assume that the quiver $Q$ is acyclic. For any vertex $t$ of $\sT_n$ and any $1\leq i\leq n$, we have
\[
X^\Z_{T_i(t)}=X^\Z_i(t).
\]
Moreover, the map taking an object $M$ to $X^\Z_M$ induces a bijection from the set of isomorphism classes of coefficient-free rigid objects of $\cD$ to the set of cluster monomials of $\cA^\Z$.
\end{Thm}
		
	\subsection{Serre polynomials}\label{sec:Serre_polynomial}
			Let $\sch/\C$ be the category of separated schemes of finite type over $\C$. For an object $X$ of $\sch/\C$, its Serre polynomial (or \emph{$E$-polynomial}) is a certain element in $\Z[q^\Hf]$ where $q^\Hf$ is an indeterminate. We refer the reader to the appendix of \cite{HauselVillegas08} for its definition.
			
			\begin{Rem} \label{Rem:pure_Serre_polynomial}
			If a scheme $X$ of $\sch/\C$ is smooth projective, then we have
			\[
			E(X)=\sum_i (-1)^i \dim H^i(X,\C) q^{\Hf i}.
			\]
			\end{Rem}
			
%\input{SerrePolynomial.tex}

% Main theorem on Serre polynomials

			A \emph{spreading out }of $X$ is given by a finitely generated $\Z$-algebra $R\subset \C$ and a separated scheme $X_R$ of finite type over $R$ \st the scheme $X$ is isomorphic to $X_R\otimes_R \C$.
			
\begin{Def}			
			An object $X$ of $\sch/\C$ is \emph{polynomial count }with \emph{counting polynomial} $P_X=P_X(t)$ in $\Z[t]$ if $X$ has a spreading out $(R,X_R)$ \st for each homomorphism $\phi:R\ra\sK$ to a finite field $\sK$, we have
			\[
			\sharp((X_R\otimes_R \sK)(\sL))=P_X(|\sL|)
			\]
			for any finite field extension $\sL\supset \sK$.
\end{Def}

\begin{Rem}
			In this case, the polynomial $P_X$ does not depend on the choice of $(R,X_R)$.
\end{Rem}

			\begin{Thm}[Katz, appendix to \cite{HauselVillegas08}]\label{thm:Katz_thm}
If an object $X$ of $\sch/\C$ is polynomial count with counting polynomial $P_X=P_X(t)$ in $\Z[t]$, then we have $E(X)=P_X(q)$ in $\Z[q^\Hf]$.
\end{Thm}
			
			\begin{Cor}\label{cor:Katz_thm}
Let $X$ be a separated scheme of finite type over $\Z$. If there exists a polynomial $P_X=P_X(t)$ in $\Z[t]$, such that we have
\[
\sharp X(\sK)=P_X(|\sK|)
\]
for every finite field $\sK$, then $E(X_\C)$ equals $P_X(q)$ in $\Z[q^\Hf]$. Here the scheme $X_\C$ is the extension of $X$ to $\C$.
\end{Cor}

By combining Remark \ref{Rem:pure_Serre_polynomial} and Theorem \ref{thm:Katz_thm}, we obtain the following result.

\begin{Cor}\label{cor:vanish_cohomology}
If an object $X$ of $\sch/\C$ is polynomial count and smooth projective. Then we have 

(1) The cohomology of $X$ is concentrated in even degrees.

(2) The Serre polynomial $E(X)$ is given by
\[
E(X)= \sum_{j} \dim H^{2j}(X,\C) q^{j}.
\]
\end{Cor}

We conclude this section with the following useful lemma.
\begin{Lem}\label{lem:q_counting_polynomial}
Let $X$ be a separated scheme of finite type defined over $\Z$. If there exists an element $P_X=P_X(v)$ of $\Z[v^\pm]$, \st we have $\sharp X(\sK)=P_X(|\sK|^{\Hf})$ for all sufficiently large finite fields $\sK$, then $P_X$ belongs to $\Z[v^2]$.
\end{Lem}
\begin{proof}
Write 
\[
P_X=f+vg,
\]
where $f$ and $g$ are Laurent polynomials in $v$ whose odd degree coefficients vanish. Suppose that the polynomial $g$ is nonzero. Then there is a (large) prime number $p$ \st $g(p^\Hf)$ does not vanish and $P_X(p^\Hf)$ equals $\sharp(X(\Fp))$. In this situation, the number $f(p^\Hf)+p^\Hf g(p^\Hf)$ is not rational, which contradicts the fact that $P_X(p^\Hf)$ is an integer. Therefore, $g$ must vanish and $P_X$ must be a Laurent polynomial in $v^2$.

We claim that $P_X$ is in fact a polynomial in $v^2$. For this, assume that we have
\[
P_X=(v^2)^{-m}h+l,
\]
where $m$ is strictly positive, $h$ is a nonzero even polynomial in $v$ of degree at most $2m-2$ and $l$ an even polynomial in $v$. Then there exists a (large) prime number $p$ \st we have $0<|h(p^\Hf)|\leq p^m-1$. This implies the inequality
\[
0<|p^{-m}h(p^\Hf)|\leq 1-p^{-m}<1,
\]
which contradicts the fact that $P_X(p^\Hf)$ and $l(p^\Hf)$ are integers.
\end{proof}

\section{Main results}\label{sec:main}

In this section we present the main results of the paper. We always consider skew-symmetric quantum cluster algebras and assume that the matrix $D$ in \eqref{eq:BZ_compatible} is the identity matrix. For convenience, we fix the following conventions: We
\begin{itemize}
\item always consider right modules,
\item denote the complex field by $k$,
\item drop the subscript $kQ$ when expressing homological data in $\mod kQ$, and
\item often denote objects and their classes in Grothendieck groups by the same symbols.
\end{itemize}

\subsection{Quantum \texorpdfstring{$F$}{F}-polynomials} \label{sec:F_M}

Let $n\geq 1$ be an integer and $Q$ an acyclic quiver with vertex set $\{1,\ldots,n\}$. Let $B$ be the associated antisymmetric matrix, \ie its entry in position $(i,j)$ is
\[
b_{ij}=\sharp\{\mathrm{arrows}\, i\ra j\}-\sharp\{\mathrm{arrows}\, j\ra i\}.
\]
The \emph{Euler form} $\langle\quad,\quad\rangle$ on $\Kz (\mod kQ)$ is defined by
\[
\langle \mU,\mV\rangle=\dim\Hom(\mU,\mV)-\dim\Ext^1(\mU,\mV)
\]
for $\mU$, $\mV$ in $\mod kQ$. If we choose the basis of $\Kz(\mod kQ)$ to be the set formed by the classes of the simples $S_i$, $1\leq i\leq n$, then $B$ is the matrix of the \emph{antisymmetrized Euler form} on $\Kz (\mod kQ)$, defined by 
\[
\langle \mU,\mV \rangle_a=\langle \mU,\mV\rangle -\langle \mU,\mV\rangle
\]
for $\mU$, $\mV$ in $\mod kQ$.

Denote by $\cC_Q$ the cluster category of $Q$, \cf \cite{BuanMarshReinekeReitenTodorov06}. Recall that there is a canonical functor $\iota:\mod kQ\ra \cC_Q$ and define
\begin{align*}
T_i=&\iota (e_i kQ),\quad 1\leq i\leq n,\\
T=&\oplus_{i=1}^n T_i.
\end{align*}
As shown in \cite{BuanMarshReinekeReitenTodorov06}, the endomorphism algebra of $T$ is canonically isomorphic to $kQ$.

Let $\sT_n$ be the $n$-regular tree with root $t_0$. For each vertex $t$ of $\sT_n$ and each $1\leq i\leq n$, we have an indecomposable rigid object $T_i(t)\in\cC_Q$ \st we have
\begin{itemize}
\item $T_i(t_0)=T_i$, $1\leq i\leq n$, and
\item if two vertices $t$ and $t'$ are linked by an edge labeled $l$, then the object $T(t')$ is obtained from $T(t)$ by the mutation at $l$.
\end{itemize}

The \emph{quantum torus }$\cR(Q)$ is the $\Rq$-algebra generated by the symbols $y^e$, $e\in \Kz(\mod kQ)$, subject to the relations
\[
y^e y^f=q^{\Hf \langle e, f\rangle _a}y^{e+f}
\]
for $e,f\in\Kz(\mod kQ)$.

Recall that we assume $Q$ to be acyclic and $E(\ )$ denotes the Serre polynomial as in Section \ref{sec:Serre_polynomial}.
\begin{Def}\label{def:F_M}
For any rigid object $M$ in $\cC_Q$, we associate to $M$ the following element in $\cR(Q)$:
\begin{align}\label{eq:F_M}
F_M=\sum_e E(\Gr_e(\Ext^1_{\cC_Q}(T,M))) q^{\Hf \langle e,e\rangle} y^e.
\end{align}
Here $\Ext^1_{\cC_Q}(T,M)$ is considered as a right $kQ$-module and we let $\Gr_e(\Ext^1_{\cC_Q}(T,M))$ denote its submodule Grassmannian whose $k$-points are the submodules with class $e$ in $\Kz(\mod kQ)$.
\end{Def}
\begin{Thm}[Quantum $F$-polynomial]\label{thm:F_M}
For any vertex $t$ of $\sT_n$ and any $1\leq i\leq n$, we have
\[
F_{T_i(t)}=F_i(t),
\]
where $F_i(t)$ is the quantum $F$-polynomial associated to $t$ and $i$ as defined in \cite{Tran09}.
\end{Thm}
We postpone the proof to Section \ref{sec:proof_F_M}.

\subsection{Quantum cluster monomials}\label{sec:X_M}

\subsubsection*{Quantum cluster algebras}

Let $m\geq n$ be two positive integers and $\tQ$ an ice quiver without loops or $2$-cycles and with vertex set $\{1,\ldots,m\}$, in which $n+1,\ldots,m$ are the frozen vertices. Let $Q$ denote the full subquiver on the vertices $1,\ldots,n$. Let $\tB$ be the $m\times n$ matrix associated to the ice quiver $\tQ$, \ie its entry in position $(i,j)$ is
\[
b_{ij}=\sharp\{\mathrm{arrows}\, i\ra j\}-\sharp\{\mathrm{arrows}\, j\ra i\}
\]
for $1\leq i\leq m$, $1\leq j\leq n$. 

Further assume that there exists some antisymmetric $m\times m$ integer matrix $\Lambda$ \st
\begin{align}\label{eq:simply_laced_compatible}
\Lambda(-\tB)=\begin{bmatrix}I_n\\0 \end{bmatrix},
\end{align}
where $I_n$ is the identity matrix of size $n\times n$. Thus, the pair $(\Lambda, \tB)$ is unitally compatible. Consequently the matrix $\tB$ is of full rank.

As described in Section \ref{sec:specialization}, to the (initial) compatible pair $(\Lambda,\tB)$, we associate the commutative cluster algebra $\cA^\Z$, the quantum cluster algebra $\cA^q$, and for each finite field $\sK$ the specialized quantum cluster algebra $\cA^\sK$. 

\subsubsection*{Triangulated categories and bilinear forms}

Let $\tW$ be a generic potential on $\tQ$. Denote the Ginzburg algebra $\Gamma(\tQ,\tW)$ by $\Gamma$. We have the triangulated categories $\per\Gamma$ and $\Dfd\Gamma$. We use the notation introduced in Section \ref{sec:cluster_category}. 

\begin{Def}[Bilinear form]\label{def:lambda}
We define a skew-symmetric bilinear form $\lambda$ on $\Kz(\per\Gamma)$ by
\[
\lambda(e_i\Gamma,e_j\Gamma)=\Lambda_{ij},\quad 1\leq i,j\leq m.
\]
\end{Def}
The following lemma follows from \eqref{eq:simply_laced_compatible} and the definitions of $\lambda$ and $\phi$.
\begin{Lem}\label{lem:lambda(P,S)}
For $1\leq i\leq m$ and $1\leq j\leq n$, we have
\[
\lambda(e_i\Gamma,\phi(S_j))=\delta_{ij},
\]
and for $1\leq i,j\leq n$, we have
\[
\lambda(\phi(S_i),\phi(S_j))=\langle S_i,S_j\rangle_a.
\]
\end{Lem}

From now on, we assume that $Q$ is acyclic. Then we have a canonical isomorphism
\[
H^0\Gamma/(e_i|n+1\leq i\leq m) \iso kQ.
\]
Accordingly there are natural inclusion functors
\[
\mod kQ\embed \mod H^0\Gamma \embed \Dfd\Gamma \embed \per\Gamma,
\]
from which we obtain maps between the Grothendieck groups
\[
\Kz(\mod kQ)\embed \Kz(\mod H^0\Gamma) \iso \Kz(\Dfd\Gamma) \ra \Kz(\per\Gamma)
\]
whose composition is denoted by $\phi$ as in Section \ref{sec:cluster_category}.

\begin{Rem}
For two finite dimensional $kQ$-modules $U$ and $V$, the $3$-Calabi-Yau property of $\Dfd\Gamma$ implies the identity
\[
\langle U,V\rangle_a=\langle U,V\rangle-\langle V,U\rangle=\langle U,V\rangle_{\Dfd\Gamma}.
\]
\end{Rem}

\subsubsection*{Refined CC-formula} 

Recall that the quantum cluster algebra $\cA^q$ is a subalgebra of the quantum torus $\cT$, where $\cT$ is the $\Rq$-algebra generated by the symbols $X^g$, $g\in \Kz(\per\Gamma)$, \st we have
\begin{itemize}
\item $X^{[e_i\Gamma]}=X_i(t_0)$ is the $i$-th initial $X$-variable in $\cA^q$, $1\leq i\leq m$, and
\item $X^g X^h =q^{\Hf \lambda(g,h)} X^{g+h}$, $g,h\in \Kz(\per\Gamma)$.
\end{itemize}

Recall that an object $M$ in $\cC$ is called coefficient-free if
\begin{enumerate}
\item  the object $M$ does not contain a direct summand $T_i$, $i>n$, and
\item  the space $\Ext^1_{\cC}(T_i,M)$ vanishes for $i>n$.
\end{enumerate}
For a coefficient-free object $M\in\cC$, the space $\Ext^1_{\cC}(T,M)$ is a right $H^0\Gamma$-module whose support is  concentrated on $Q$. Thus, it can be viewed as a $kQ$-module.

Recall that the quiver $Q$ is acyclic and that $E(\ )$ denotes the Serre polynomial as in Section \ref{sec:Serre_polynomial}.
\begin{Def}[Refined CC-formula]\label{def:X_M}
For any coefficient-free rigid object $M$ of $\cD$, we denote by $m$ the class of $\Ext^1_{\cC}(T,M)$ in $\Kz(\mod kQ)$, and associate to $M$ the following element in $\cT$:
\begin{align}\label{eq:X_M}
X_M=\sum_e E(\Gr_e (\Ext^1_{\cC}(T,M))) q^{-\Hf \langle e,m-e \rangle} X^{\ind_T(M)-\phi(e)}.
\end{align}
\end{Def}

\begin{Thm}[Main theorem]\label{thm:X_M}
For any vertex $t$ of $\sT_n$ and any $1\leq i\leq n$, we have
\[
X_{T_i(t)}=X_i(t).
\]
Moreover, the map taking an object $M$ to $X_M$ induces a bijection from the set of isomorphism classes of coefficient-free rigid objects of $\cD$ to the set of quantum cluster monomials of $\cA^q$.
\end{Thm}

\begin{Thm}[Geometry of quiver Grassmannians]\label{thm:Grassmannian_geometry}\footnote{This theorem first appeared in \cite{CalderoReineke08}. But as pointed out by \cite{Nakajima09}, the proof there contains a gap. Nevertheless, the assertion holds true, and we shall give a proof using quantum cluster algebras.}
Assume that the quiver $Q$ is acyclic. For any rigid $kQ$-module $V$ and each class $e$ of $\Kz(\mod kQ)$, the complex variety $\Gr_e(V)$ has the following properties:
\begin{enumerate}
\item Its cohomology is concentrated in even degrees.
\item Its Serre polynomial counts the number of its rational points over any finite field $\sK$.
\end{enumerate}
\end{Thm}

\begin{Rem}
As we shall see, our proof of Theorem \ref{thm:Grassmannian_geometry} uses combinatorial data of quantum cluster algebras (specialization maps). On the other hand, a purely geometric proof of part (1) is also available, due to Hiraku Nakajima, where one uses Fourier-Sato-Deligne transform and the Decomposition theorem. More details of this geometric approach will appear in a new version of \cite{Nakajima09}.
\end{Rem}

We postpone the proofs of Theorem \ref{thm:X_M} and Theorem \ref{thm:Grassmannian_geometry} to Section \ref{sec:proof_acyclic} and Section \ref{sec:proof_general_coefficient}.

\begin{Cor}[Positivity]
The cluster monomials in $\cA^\Z$ have non-negative coefficients in their expansions as Laurent polynomials in the initial variables.
\end{Cor}
\begin{proof}
The map taking an object $M$ to
\[
X_M^{\Z}=\sum_e \chi(\Gr_e\Ext^1_\cC(T,M))X^{\ind_T M-\phi(e)}
\]
induces a bijection from the set of isomorphism classes of coefficient-free rigid objects in $\cD$ to the set of the cluster monomials in $\cA^\Z$, \cf Theorem \ref{thm:CC_formula}. In addition, for any coefficient-free rigid object $M$ the $kQ$-module $\Ext^1_{\cC}(T,M)$ is rigid. Therefore the claim follows from Part (1) of Theorem \ref{thm:Grassmannian_geometry}.
\end{proof}

\subsection{Bilinear form}
In this section, let $\lambda$ be the bilinear form as in Definition \ref{def:lambda}. The quiver $Q$ is not assumed to be acyclic.

\subsubsection*{\texorpdfstring{$\Lambda$}{Lambda}-matrices}

Let $t$ and $t'$ be vertices of $\sT_n$ linked by an edge labeled $k$. Let $\Lambda(t)$ and $\Lambda(t')$ be the $\Lambda$-matrices in $t$ and $t'$ respectively.
\begin{Lem}\label{lem:basis_change}
Let $X\in\cC$ be rigid. Then for some sign $\epsilon$, we have
\[
[\ind_{T(t')} X:T_i(t')]=\sum_{1\leq j\leq m} (E_{\epsilon})_{ij}[\ind_{T(t)} X:T_j(t)],
\]
where $E_{\epsilon}$ is the matrix associated with $\tB(t)$ and $k$ in equation \eqref{eq:E_epsilon}.
\end{Lem}

\begin{proof}\footnote{When $(\tQ,\tW)$ is Jacobi-finite, the claim also follows from Section $3$ of \cite{DehyKeller08}.}
Lift $T(t)$ and $T(t')$ to $\pi^{-1}T(t)$ and $\pi^{-1}T(t')$ in $\cF$. By \cite[Theorem 2.18]{Plamondon10a} the objects $\pi^{-1}T(t)$ and $\pi^{-1}T(t')$ are linked by a mutation in the sense of \cite{KellerYang09}. Therefore, when we choose the classes $[\pi^{-1}T_i(t)]$, $1\leq i\leq m$, respectively $[\pi^{-1}T_i(t')]$, $1\leq i\leq m$, as a basis of $\Kz(\per\Gamma)$, the corresponding change of coordinates matrix equals $E_\epsilon$ for some sign $\epsilon$.
\end{proof}

\begin{Prop}[Bilinear form]\label{prop:matrix_representation}
For any vertex $t$ of $\sT_n$, the matrix $\Lambda(t)$ is the matrix of $\lambda$ in the basis formed by the classes $[\pi^{-1}T_i(t)]$, $1\leq i\leq m$.
\end{Prop}
\begin{proof}
Recall that if $t$ and $t'$ are neighboring vertices of $\sT_n$, for any sign $\epsilon$, we have $\Lambda(t')=E_{\epsilon}^T \Lambda(t) E_{\epsilon}$. Starting from $t=t_0$ the root of $\sT_n$, we inductively deduce the claim from Lemma \ref{lem:basis_change}.
\end{proof}

\subsubsection*{Euler form}
For an object $P$ of $\per\Gamma$ and an object $S$ of $\Dfd\Gamma$, the spaces $\Hom_{\per\Gamma}(P,S[i])$ are finite dimensional and vanish for all but finitely many integers $i$. Therefore, the Euler form given by
\[
\sum_{i\in\Z}(-1)^i\dim\Hom_{\per\Gamma}(P,S[i])
\]
is a well defined bilinear form, which we denote by
\[
\langle\ ,\ \rangle_{\per\Gamma}:\Kz (\per\Gamma) \times \Kz(\Dfd\Gamma) \Ra \Z.
\]
We define the \emph{support} of an object $S$ of $\Dfd\Gamma$ to be the set of vertices $i$ \st the space $\Hom_{\per\Gamma}(e_i\Gamma,S)$ does not vanish.
\begin{Prop} \label{prop:bilinear_form}
For an object $P$ of $\per\Gamma$ and an object $S$ of $\Dfd\Gamma(Q,W)$ whose support is contained in $\{1,\ldots,n\}$, we have
\[
\lambda([P],\phi([S]))=\langle [P],[S]\rangle_{\per\Gamma}.
\]
\end{Prop}
\begin{proof}
The claim follows from Lemma \ref{lem:lambda(P,S)}.
\end{proof}

\section{Examples}

\subsection{\texorpdfstring{$A_2$}{A2} with coefficients}

Let $\tQ$ be the ice quiver defined as
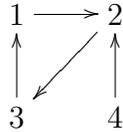
\[
\xymatrix{
 1 \ar[r]& 2  \ar[dl]\\
 3 \ar[u]& 4  \ar[u]}
\]
where the frozen vertices are $3$ and $4$. Let $\tB$ be the $4\times 2$ matrix associated to the ice quiver $\tQ$. Let $\tB^\circ$ be the $4\times 4$ matrix associated to the ordinary quiver $\tQ^\circ$ underlying $\tQ$. It is invertible and we define a matrix $\Lambda=(\tB^\circ)^{-T}$. Then the pair $(\Lambda, \tB)$ is compatible. Associate to it a quantum cluster algebra $\cA^q$.

Denote by $F$ the functor $\Ext^1_\cC(T,\ )$. We associate to the $2$-regular tree $\sT_2$ with root $t_0$
\[
\xymatrix{
\cdots t_{-1}\ar@{-}[r]^2 & t_0\ar@{-}[r]^1 &t_1\ar@{-}[r]^2&t_2\ar@{-}[r]^1&t_3\ar@{-}[r]^2&t_4\ar@{-}[r]^1&t_5\ar@{-}[r]^2&t_6\cdots
}
\]

the following quantum cluster variables

\begin{table}[h!]
\begin{center}
\begin{tabular}{ |c |c|l|c|}
\hline
$t$& $i$&$X_i(t)$& $\dimv F(T_i(t))$\\
\hline
$t_1 $&$1$&\small X(-1,0,1,0)+X(-1,1,0,0) & (1,0,0,0)\\
\hline
$t_2 $&$2$& \small { X(-1,-1,1,0)+X(-1,0,0,0)+X(0,-1,0,1)}&(1,1,0,0)\\
\hline
$t_3 $&$1$& \small X(0,-1,1,0)+X(1,-1,0,1) & (0,1,0,0)\\
\hline
$t_4 $&$2$&\small X(1,0,0,0) & (0,0,0,0)\\
\hline
$t_5 $&$1$&\small X(0,1,0,0) & (0,0,0,0)\\
\hline
\end{tabular}
\caption{Quantum cluster variables in case $A_2$ with coefficients.}
\end{center}
\end{table}

One easily checks that Theorem \ref{thm:X_M} holds for these variables, as every submodule Grassmannian appearing in our case has at most one $\C$-point.

\subsection{A regular rigid module}

Let $Q$ be of extended Dynkin type $\tilde{E_6}$ with the following orientation
\begin{align*}
\begin{CD}
  @.          @.     5 			@.          @.  \\
@.      @.        @AAA                 @.\\
 @.    		@. 	 4 	@. 			   @.  	\\
@.	    @. 			  @AAA					@.\\
3 @<<<   2    @<<< 	  1   	   			@>>> 			6 	@>>> 	7.
\end{CD}
\end{align*}

Let $\tQ$ be the ice quiver obtained from $Q$ by adding the frozen vertices $\{8,9,\ldots,14\}$ and one arrow $i+7\ra i$ for each $1\leq i\leq 7$. Let $\tB$ be the $14\times 7$ matrix associated to the ice quiver $\tQ$. Let $\tB^\circ$ be the $14\times 14$ matrix associated to the ordinary quiver $\tQ^\circ$ underlying $\tQ$. It is invertible and we define a matrix $\Lambda=(\tB^\circ)^{-T}$. Then the pair $(\Lambda, \tB)$ is compatible. Associate to it a quantum cluster algebra $\cA^q$.

Consider the object $T(t)$ defined as $\mu_{1}\mu_{ 6}\mu_{ 4}\mu_{ 2}\mu_{ 1}\mu_{ 7}\mu_{ 5}\mu_{ 3}T(t_0)$. Denote the class of $\Ext^1_{\cC}(T,T_1(t))$ in $K_0\mod kQ$ by $m$. A computer program gives the results in Table \ref{tab:Regular}.

%%%%%%%%%%
%%what's the better way to ajust the position of this table?
%%%
\begin{table}[htbp]\label{tab:Regular}
\begin{center}
\begin{tabular}{  | p{12cm} |}
    \hline
         $\dimv m$  \\ \hline(2,1,1,1,1,1,1)\\ \hline
    Extended $g$-vector $\tg_1(t)$ \\ \hline $-2e_1+e_2+e_4+e_6$\\ \hline
    Commutative $F$-polynomial $F_1^\Z(t)$\\ \hline 
    $1 $+$ 2 {y_1} $+$ {y_1}^2 $+$ ({y_1} $+$ {y_1}^2) {y_2} $+$ (({y_1} $+$ {y_1}^2) {y_2}) {y_3} $+$ ({y_1} $+$ {y_1}^2 $+$ {y_1}^2 {y_2} $+\\
    $ {y_1}^2 {y_2} {y_3}) {y_4} $+$ (({y_1} $+$ {y_1}^2 $+$ {y_1}^2 {y_2} $+$ {y_1}^2 {y_2} {y_3}) {y_4}) {y_5} $+$ ({y_1} $+$ {y_1}^2 $+$ {y_1}^2 {y_2} $\\
    +$ {y_1}^2 {y_2} {y_3} $+$ ({y_1}^2 $+$ {y_1}^2 {y_2} $+$ {y_1}^2 {y_2} {y_3}) {y_4} $+$ (({y_1}^2 $+$ {y_1}^2 {y_2} $+$ {y_1}^2 {y_2} {y_3}) {y_4}) {y_5}) {y_6} $\\
    +$ (({y_1} $+$ {y_1}^2 $+$ {y_1}^2 {y_2} $+$ {y_1}^2 {y_2} {y_3} $+$ ({y_1}^2 $+$ {y_1}^2 {y_2} $+$ {y_1}^2 {y_2} {y_3}) {y_4} $+$ (({y_1}^2 $+$ {y_1}^2 {y_2} $\\
    +$ {y_1}^2 {y_2} {y_3}) {y_4}) {y_5}) {y_6}) {y_7}$ \\
\hline
    Quantum cluster variable $X_1(t)$ \\ \hline
\small X(-2,-1,0,-1,0,-1,0,2,0,0,0,0,0,0)+$(q^{-\Hf}+q^{\Hf})$X(-2,0,0,0,0,0,0,1,0,0,0,0,0,0)\\
\small+X(-2,1,0,1,0,1,0,0,0,0,0,0,0,0)+X(-1,-1,-1,-1,0,-1,0,2,1,0,0,0,0,0)\\
\small+X(-1,-1,0,-1,-1,-1,0,2,0,0,1,0,0,0)+X(-1,-1,0,-1,0,-1,-1,2,0,0,0,0,1,0)\\
\small+X(-1,-1,0,-1,0,0,-1,2,0,0,0,0,1,1)+X(-1,-1,0,0,-1,-1,0,2,0,0,1,1,0,0)\\
\small+X(-1,0,-1,-1,0,-1,0,2,1,1,0,0,0,0)+X(-1,0,-1,0,0,0,0,1,1,0,0,0,0,0)\\
\small+X(-1,0,0,0,-1,0,0,1,0,0,1,0,0,0)+X(-1,0,0,0,0,0,-1,1,0,0,0,0,1,0)\\
\small+X(-1,0,0,0,0,1,-1,1,0,0,0,0,1,1)+X(-1,0,0,1,-1,0,0,1,0,0,1,1,0,0)\\
\small+X(-1,1,-1,0,0,0,0,1,1,1,0,0,0,0)+X(0,-1,-1,-1,-1,-1,0,2,1,0,1,0,0,0)\\
\small+X(0,-1,-1,-1,0,-1,-1,2,1,0,0,0,1,0)+X(0,-1,-1,-1,0,0,-1,2,1,0,0,0,1,1)\\
\small+X(0,-1,-1,0,-1,-1,0,2,1,0,1,1,0,0)+X(0,-1,0,-1,-1,-1,-1,2,0,0,1,0,1,0)\\
\small+X(0,-1,0,-1,-1,0,-1,2,0,0,1,0,1,1)+X(0,-1,0,0,-1,-1,-1,2,0,0,1,1,1,0)\\
\small+X(0,-1,0,0,-1,0,-1,2,0,0,1,1,1,1)+X(0,0,-1,-1,-1,-1,0,2,1,1,1,0,0,0)\\
\small+X(0,0,-1,-1,0,-1,-1,2,1,1,0,0,1,0)+X(0,0,-1,-1,0,0,-1,2,1,1,0,0,1,1)\\
\small+X(0,0,-1,0,-1,-1,0,2,1,1,1,1,0,0)+X(1,-1,-1,-1,-1,-1,-1,2,1,0,1,0,1,0)\\
\small+X(1,-1,-1,-1,-1,0,-1,2,1,0,1,0,1,1)+X(1,-1,-1,0,-1,-1,-1,2,1,0,1,1,1,0)\\
%\hline
%\end{tabular}
%\end{center}
%\end{table}
%\begin{table}[htbp]
%\begin{center}
%\begin{tabular}{  | p{12cm} |}
%    \hline
\small+X(1,-1,-1,0,-1,0,-1,2,1,0,1,1,1,1)+X(1,0,-1,-1,-1,-1,-1,2,1,1,1,0,1,0)\\
\small+X(1,0,-1,-1,-1,0,-1,2,1,1,1,0,1,1)+X(1,0,-1,0,-1,-1,-1,2,1,1,1,1,1,0)\\
\small+X(1,0,-1,0,-1,0,-1,2,1,1,1,1,1,1)\\
\hline
\end{tabular}
\caption{Cluster data corresponding to a regular rigid module with a non-trivial submodule Grassmannian.}
\end{center}
\end{table}

We claim that $X_{T_1(t)}=X_1(t)$. To see this, recall that we are working with right $\C Q$-modules, so the only submodule Grassmannian of $\Ext^1_{\cC}(T,T_1(t))$ containing more than one $\C$-point is $\Gr_{[S_1]}\Ext^1_{\cC}(T,T_1(t))$. It contributes the only nonzero coefficient other than $1$ in $X_{T_1(t)}$, which is 
\begin{align*}
E(\Gr_{[S_1]}\Ext^1_{\cC}(T,T_1(t)))q^{-\Hf \langle[S_1],[\Ext^1_{\cC}(T,T_1(t))]-[S_1]\rangle}=&E(\P ^1)q^{-\Hf }\\
=&q^{-\Hf}+q^{\Hf}.
\end{align*}
This coefficient agrees with the only nonzero coefficient other than $1$ in $X_1(t)$.

\subsection{A generic non-rigid module}

Submodule Grassmannians do not always have counting polynomials, even if they happen to be smooth. Let us consider an ice quiver $\tQ$
\[
\xymatrix{
1\ar@<.9ex>[rr]\ar@<.3ex>[rr]\ar@<-.3ex>[rr]\ar@<-.9ex>[rr] & & 2\ar[d]\\
3\ar[u] & & 4}
\]
whose principal part $Q$ is the full subquiver supported on the vertex set $\{1,\ 2\}$.

As in Example 3.6 of \cite{DerksenWeymanZelevinsky09}, let $M$ be a generic representation of $Q$ with dimension vector $(3,4)$. As shown in [loc. cit.], the submodule Grassmannian consisting of subrepresentations with dimension vector $(1,1)$ is smooth. If this variety had a counting polynomial, applying Corollary \ref{cor:vanish_cohomology} we would obtain the vanishing of its odd cohomology and thereby the positivity of its Euler characteristic. On the other hand, its Euler characteristic has been computed in [loc. cit.] and turns out to be $-4$. This contraditiction shows that this smooth submodule Grassmannian does not have a counting polynomial.

\section{Proof for acyclic coefficients}\label{sec:proof_acyclic}

For convenience, for a fixed finite field $\sK$ and any scheme $X$ of finite type over $\Z$, we often denote the number of $\sK$-points of $X$ by $\sharp(X)$.

\subsection{Theorems for acyclic coefficients} \label{sec:acyclic_main_thm} If the ice quiver $\tQ$ is acyclic, then the presentable cluster category $\cD$ agrees with $\cC$ and the usual cluster category of $\tQ$ as in \cite{BuanMarshReinekeReitenTodorov06}. Our main aim in this section is to prove the following two theorems.

\begin{Thm}[Main theorem for acyclic coefficients]\label{thm:acyclic_main_thm}
Assume that the ice quiver $\tQ$ is acyclic. For any vertex $t$ of $\sT_n$ and any $1\leq i\leq n$, we have
\[
X_{T_i(t)}=X_i(t).
\]
Moreover, the map taking an object $M$ to $X_M$ induces a bijection from the set of isomorphism classes of coefficient-free rigid objects of $\cC$ to the set of quantum cluster monomials of $\cA^q$.
\end{Thm}

\begin{Thm*}[\ref{thm:Grassmannian_geometry}]
Assume that the quiver $Q$ is acyclic. For any rigid $kQ$-module $V$ and each class $e$ of $\Kz(\mod kQ)$, the complex variety $\Gr_e(V)$ has the following properties:
\begin{enumerate}
\item Its cohomology is concentrated in even degrees.
\item Its Serre polynomial counts the number of its rational points over any finite field $\sK$.
\end{enumerate}
\end{Thm*}

\begin{proof}[Proofs of Theorem \ref{thm:acyclic_main_thm} and \ref{thm:Grassmannian_geometry}] Our strategy is inspired by the counting property of Serre polynomials in Corollary \ref{cor:Katz_thm}. Instead of attacking Theorem \ref{thm:acyclic_main_thm} directly, we look for counting polynomials with the help of specialization maps on quantum cluster algebras. The proof consists of the following steps.

(i) For each finite field $\sK$, consider the specialized quantum cluster algebra $\cA^{\sK}$. Recall that there is a specialization map $ev_\sK$ from $\cA^q$ to $\cA^\sK$. We will show in Theorem \ref{thm:spec_X_M} that the specialized version of Theorem \ref{thm:acyclic_main_thm} is true. Namely, to each coefficient-free rigid object $M$, associate an element
\[
X^\sK_M=\sum_e \sharp(\Gr_e \Ext^1_\cC(T,M))|\sK|^{-\Hf \langle e,m-e\rangle}X^{\ind_T M-\phi(e)}.
\]
Notice that the $kQ$-module $\Ext^1_\cC(T,M)$ is rigid so that the Grassmannian $\Gr_e(\Ext^1_\cC(T,M))$ is canonically defined over $\Z$, \cf Section \ref{sec:int_rigid_mod}.

We will show that for any vertex $t$ of $\sT_n$ and any $1\leq i\leq n$, we have
\[
X^\sK_{T_i(t)}=X_i^\sK(t),
\]
and that the map taking an object $M$ to $X^\sK_M$ induces a bijection from the set of isomorphism classes of coefficient-free rigid objects of $\cC$ to the set of quantum cluster monomials of $\cA^\sK$.

(ii) We proceed to find counting polynomials for submodule Grassmannians. By Theorem \ref{thm:CC_formula}, every cluster monomial takes the form
\[
X_M^\Z=\sum_e \chi(\Gr_e\Ext^1_\cC(T,M)) X^{\ind_T M-\phi(e)}
\]
for some coefficient-free rigid object $M$. Recall from \ref{prop:spec_bijective} that the specialization map $ev_1:\cA^q\ra\cA^\Z$ induces a support-preserving bijection between the sets of (quantum) cluster monomials. Since the map $\phi$ is injective, we can write the quantum cluster monomial corresponding to $X^\Z_M$ in the form
\[
ev_1^{-1}X^\Z_M=\sum_e P_e X^{\ind_T M-\phi(e)},
\]
where the coefficients $P_e$ belong to $\Z[q^{\pm\Hf}]$ and have subtraction-free rational expressions in $q^\Hf$. The specialization map $ev_{\sK}:\cA^q\ra\cA^{\sK}$ sends  $ev_1^{-1}X^\Z_M$ to the quantum cluster monomial $ev_\sK ev_1^{-1} (X^\Z_M)$ of $\cA^\sK$ which has the unique minimal degree $\ind_T M$ in the sense of Section \ref{sec:specialization}. Proposition \ref{prop:minimal_degree} claims that different quantum cluster monomials of $\cA^\sK$ have different minimal degrees. So the quantum cluster monomials $ev_\sK ev_1^{-1} (X^\Z_M)$ and $X^\sK_M$ must be the same since they have the same minimal degree, \ie we have
\begin{align*}
\begin{split}
&\sum_e ev_\sK(P_e) X^{\ind_T M-\phi(e)}\\
&\qquad=\sum_e \sharp(\Gr_e\Ext^1_\cC(T,M))|\sK|^{-\Hf \langle e,m-e\rangle} X^{\ind_T M-\phi(e)}.
\end{split}
\end{align*}
Thus, since the map $\phi$ is injective, we obtain for each finite field $\sK$ the identity
\begin{align*}
P_e(|\sK|^\Hf)|\sK|^{\Hf \langle e,m-e\rangle}=\sharp(\Gr_e\Ext^1_\cC(T,M)).
\end{align*}
Lemma \ref{lem:q_counting_polynomial} implies that $P_e q^{\Hf \langle e,m-e\rangle}$ belongs to $\Z[q]$. We denote this polynomial by $E_e$.

(iii) From Corollary \ref{cor:Katz_thm} we deduce that the counting polynomial $E_e$ agrees with the Serre polynomial $E(\Gr_e\Ext^1_\cC(T,M))$. So Part (2) of Theorem \ref{thm:Grassmannian_geometry} and Theorem \ref{thm:acyclic_main_thm} hold true. 

(iv) Finally, because $X=\Gr_e\Ext^1_\cC(T,M)$ is smooth projective, by Remark \ref{Rem:pure_Serre_polynomial}, we have
\[
E(X)=\sum_i (-1)^i \dim H^i(X,\C) q^{\Hf i}.
\]
Since $E(X)$ is a polynomial in $q$, the cohomology $H^i(X,\C)$ must vanish for odd degrees $i$, which shows Part (1) of Theorem \ref{thm:Grassmannian_geometry}.
\end{proof}

\subsection{\texorpdfstring{$\Lambda$}{Lambda}-matrices} 
We proceed to compute some useful $\Lambda$-matrices. Assume that the principal part $Q$ is acyclic.

\begin{Lem} \label{lem:<e,m>} For any coefficient-free rigid object $M$ in $\cD$ which contains no summands in $\add T$, and for any classes $e$, $a$, $b$ in $\Kz(\mod kQ)$, we have
\begin{align}
&\lambda(\ind_T M,-\phi(e))=-\ind_T M\cdot e= \langle e,m \rangle = \langle e,m-e \rangle + \langle e,e \rangle, \\
&\lambda(-\phi(a),-\phi(b))= \langle a,b \rangle - \langle b,a \rangle .
\end{align}
Here we let $m$ denote the class of $\Ext^1_\cC(T,M)$ in $\Kz(\mod kQ)$, $e_i$ denote the $i$th coordinate of $e$ under the basis of $\Kz(\mod kQ)$ formed by $[S_i]$, $1\leq i\leq n$, and define the product
\[
\ind_T M\cdot e=\sum_{1\leq i\leq n} [\ind_T M: T_i] e_i.
\]
\end{Lem}
\begin{proof}
The equality $-\ind_T M\cdot e=\langle e,m\rangle$ is a immediate consequence from Lemma 4.7 of \cite{Plamondon10b}. The rest of the assertions follow from Lemma \ref{lem:lambda(P,S)}.
\end{proof}

\begin{Cor} \label{cor:expand_lambda}For $\mM$, $\mN$ as in the above lemma, we have
\begin{align*}
&\Hf \lambda(\ind_T \mM -\phi(\db),\ind_T \mN -\phi(\dd))\\
&\qquad=\Hf \lambda(\ind_T \mM ,\ind_T \mN )+\Hf  \langle \dd,\dm \rangle -\Hf  \langle \db,\dn \rangle +\Hf ( \langle \db,\dd \rangle - \langle \dd,\db \rangle ).
\end{align*}
\end{Cor}

\begin{Lem} For any coefficient-free rigid object $M$ in $\cD$, and any class $e$ which is the class of some submodule of $\Ext^1_\cC(T,M)$, we have
\[
\Hf\langle e,e \rangle=-\Hf \langle e,m-e \rangle +\Hf \lambda(\ind_T M, -\phi(e))
\]
\end{Lem}
\begin{proof}
Since $\langle e,e \rangle=-\langle e,m-e \rangle +\langle e, m\rangle $, it suffices to prove
\[
\langle e, m\rangle=\lambda(\ind_T M, -\phi(e)).
\]

(i) If the object $M$ contains no summand of $\add T$, the assertion follows from Lemma \ref{lem:<e,m>}.

(ii) If the object $M$ belongs to $\add T$, we have $m=e=0$ and the assertion holds.

(iii) For an object $M$ in the form $M=M_1\oplus M_2$, where $M_1$ contains no summand of $\add T$ and $M_2$ belongs to $\add T$, since the object $M$ is rigid, we have
\begin{align*}
\sum_{1\leq i\leq n}[\ind_T M_2:T_i] \cdot m_i&=\dim\Hom(\Ext^1_\cC(T,\tau^{-1}M_2), \Ext^1_\cC(T,M_1))\\
&\qquad \leq \dim \Ext^1_\cC(M_2,M_1)=0,
\end{align*}
and consequently the identity
\[
\lambda(\ind_T M_2, -\phi(e))=-\sum_{1\leq i\leq n}[\ind_T M_2:T_i] \cdot e_i=0.
\]
Also, step (i) implies the identity
\[
\lambda(\ind_T M_1, -\phi(e))=\langle e, m\rangle.
\]
Then the claim follows.
\end{proof}

This Lemma implies the following results.
\begin{Def}\label{def:general_F_M}
For any coefficient-free rigid object $M$ of $\cD$, we associate to $M$ the following element in $\cR(Q)$:
\[
F_M=\sum_e E(\Gr_e \Ext^1(T,M))q^{\Hf\langle e,e\rangle}  y^e.
\]
\end{Def}
\begin{Prop}\label{prop:decompose_X_M}
For any coefficient-free rigid object $M$ of $\cD$, we have
\[
X_M=F_M|_{y^e\mapsto  X^{-\phi(e)}}X^{\ind_T M}.
\]
\end{Prop}

%temp2.tex
From now on, assume that the ice quiver $\tQ$ is acyclic.

\begin{Thm}
\cite[Theorem 7.5]{BuanMarshReinekeReitenTodorov06} For any pair of indecomposable rigid objects $\mN, \mM\in \cC$ \st the dimension of the space $\Ext^1_\cC(\mN,\mM)$ equals $1$, there exists an integer $1\leq k\leq n$, and two vertices $t_1$ and $t_2$ of $\sT_n$ linked by an edge labeled $k$, \st we have $\mN=T_k(t_1)$, $\mM=T_k(t_2)$.
\end{Thm}
Accordingly, there are two non-split triangles in $\cC$
\begin{align}
\mN\ra \mBp\ra \mM\ra \Sigma \mN,\label{tri:lift}\\
\mM\ra \mBn\ra \mN\ra \Sigma \mM,
\end{align}
where the second terms are
\[
\mBp= \oplus_{1\leq j\leq m} T_j(t_1)^{[-b_{jk}(t_1)]_{+}},
\]
and
\[
\mBn= \oplus_{1\leq i\leq m} T_i(t_1)^{[b_{ik}(t_1)]_{+}}.
\]
Notice that the objects $T_i(t_1)$ and $T_i(t_2)$ are the same for $1\leq i\leq m$, $i\neq k$.
 
According to \cite[Proposition 2.21]{Plamondon10a}, exactly one of the above triangles lifts to the fundamental domain $\cF$ of $\per\Gamma$. Assume that the triangle \eqref{tri:lift} lifts. Then we have
\begin{align*}
\ind_T \mBp=&\ind_T \mN+\ind_T \mM,\\
\ind_T \mBn=&\ind_T \mBp+(\ind_T \mBn-\ind_T \mBp)\\
						=&\ind_T \mBp + \sum_i b_{ik}(t_1)\ind_T T_i(t_1),\\
\ind_T \mBp=&\ind_T \mBn -  \sum_i  b_{ik}(t_1)\ind_T T_i(t_1)\\
						=&\ind_T \mBn +  \sum_i b_{ik}(t_2)\ind_T T_i(t_2).
\end{align*}

\begin{Lem} \label{lem:compute_lambda} For the objects $\mN$, $\mM$ as above \st $\ind_T \mBp$ equals $\ind_T \mN+\ind_T \mM$, we have
\begin{align*}
\Lambda(t_1)(e_k, \sum_j [-b_{jk}(t_1)]_{+}e_j)&=\lambda({\ind_T  \mN},{\ind_T \mBp})=\lambda({\ind_T \mN},{\ind_T \mM}),
\\
\Lambda(t_1)(e_k, \sum_i [b_{ik}(t_1)]_{+}e_i)&=\lambda({\ind_T \mN},{\ind_T \mBn})=\lambda({\ind_T \mN},{\ind_T \mM})-1,
\\
\Lambda(t_2)(e_k, \sum_j [b_{jk}(t_2)]_{+}e_j)&=\lambda({\ind_T \mM},{\ind_T \mBp})=\lambda({\ind_T \mM},{\ind_T \mN}),
\\
\Lambda(t_2)(e_k, \sum_i [-b_{ik}(t_2)]_{+}e_i)&=\lambda({\ind_T \mM},{\ind_T \mBn})\\
&\qquad=\lambda({\ind_T \mM},{\ind_T \mN})+1.
\end{align*}
\end{Lem}

%temp1.tex

\begin{proof}
Since the bilinear form $\lambda$ is skew-symmetric, we have
\begin{align*}
\lambda({\ind_T  \mN},{\ind_T \mBp})&=\lambda({\ind_T \mN},{\ind_T \mM}),\\
\lambda({\ind_T \mN},{\ind_T \mBn})&=\lambda({\ind_T \mN},{\ind_T \mM})+\lambda({\ind_T \mN},\sum_i b_{ik}(t_1)\ind_T T_i(t_1)),\\
\lambda({\ind_T \mM},{\ind_T \mBp})&=\lambda({\ind_T \mM},{\ind_T \mN}),\\
\lambda({\ind_T \mM},{\ind_T \mBn})&=\lambda({\ind_T \mM},{\ind_T \mN})\\
&\qquad-\lambda({\ind_T \mM},\sum_i b_{ik}(t_2)\ind_T T_i(t_2)).
\end{align*}

Applying Proposition \ref{prop:matrix_representation} we obtain
\begin{align*}
&\Lambda(t_1)(e_k, \sum_j [-b_{jk}(t_1)]_{+}e_j)=\lambda({\ind_T  \mN},{\ind_T \mBp}),\\
&\Lambda(t_1)(e_k, \sum_i [b_{ik}(t_1)]_{+}e_i)=\lambda({\ind_T \mN},{\ind_T \mBn}),\\
&\Lambda(t_2)(e_k, \sum_j [b_{jk}(t_2)]_{+}e_j)=\lambda({\ind_T \mM},{\ind_T \mBp}),\\
&\Lambda(t_2)(e_k, \sum_i [-b_{ik}(t_2)]_{+}e_i)=\lambda({\ind_T \mM},{\ind_T \mBn}),\\
&\lambda({\ind_T \mN},\sum_i b_{ik}(t_1)\ind_T T_i(t_1))=\Lambda(t_1)(e_k,b^k(t_1)),\\
&\lambda({\ind_T \mM},\sum_i b_{ik}(t_2)\ind_T T_i(t_2))=\Lambda(t_2)(e_k,b^k(t_2)).
\end{align*}
Here $\tb^k(t_1)$ and $\tb^k(t_2)$ are the $k$-th columns of $\tB(t_1)$ and $\tB(t_2)$ respectively.

Finally, recall
\[
\Lambda(t_1)(-\tB(t_1))=\Lambda(t_2)(-\tB(t_2))=\begin{bmatrix}I_n\\0 \end{bmatrix},
\]
 we obtain
\begin{align*}
&\Lambda(t_1)(e_k,b^k(t_1))=-1,\\
&\Lambda(t_2)(e_k,b^k(t_2))=-1,
\end{align*}
and thus complete the proof.
\end{proof}

\subsection{Specialized formula}\label{sec:spec_X_M}

Recall that the specialized quantum cluster algebra $\cA^\sK$ is a $\Rk$-subalgebra of the quantum torus $\cT^\sK$, where $\cT^\sK$ is the $\Rk$-algebra generated by the symbols $X^g$, $g\in \Kz(\per\Gamma)$, \st
\begin{itemize}
\item $X^{[e_i\Gamma]}=X_i(t_0)$ is the $i$-th initial $X$-variable in $\cA^\sK$, $1\leq i\leq m$, and we have
\item $X^g X^h =|\sK|^{\Hf \lambda(g,h)} X^{g+h}$, $g,h\in \Kz(\per\Gamma)$.
\end{itemize}
We fix the following conventions in this subsection: 
\begin{itemize}
\item We denote the functor $\Ext^1_{\cC}(T,\ )$ by $F$.
\item We assume that the principal part $Q$ is acyclic.
\item We shall use uppercase letters to denote $\sK Q$-modules and lowercase letters to denote their dimension vectors. We choose the the set formed by the $[S_i]$, $1\leq i\leq n$ as a basis of $\Kz(\mod \sK Q)$ and identify the dimension vectors with classes in $\Kz(\mod \sK Q)$.
\end{itemize}
We propose the following formula:
\begin{Def}[Specialized formula]\label{def:spec_X_M}
For any coefficient-free and rigid object $M\in\cD$, we denote by $m$ the class of $FM$ in $\Kz(\mod kQ)$, and associate to $M$ the following element in $\cT^\sK$:
\begin{align}\label{eq:spec_X_M}
X^\sK_M=\sum_e \sharp(\Gr_e FM) |\sK|^{-\Hf \langle e,m-e \rangle} X^{\ind_T(M)-\phi(e)}.
\end{align}
\end{Def}

In this subsection, we shall prove that the formula $X^\sK_?$ above takes direct sums to products.

\begin{Prop}[Multiplicativity]\label{prop:direct_sum}
For any coefficient-free rigid object $ \mM \oplus  \mN $ of $\cD$, we have
\begin{align}\label{eq:spec_mult}
X^\sK_ \mM \cdot X^\sK_ \mN=|\sK|^{\Hf \lambda(\ind_T \mM ,\ind_T \mN )}X^\sK_{ \mM \oplus  \mN }
\end{align}
\end{Prop}
\begin{proof}

In order to prove the assertion, we distinguish five cases:

(i) Assume that the object $\mM\oplus\mN$ contains no summands from $\add T$.

For any submodule $\mD$ of $F\mN$ and any submodule $\mB$ of $F\mM$, the fact $\Ext^1(F\mN,F\mM)=0$ implies that $\Ext^1(D,F\mM/B)$ vanishes. Denote by $\dn$ the class of $F\mN$ in $\Kz(\mod kQ)$. Then the argument of Lemma 3.8 of \cite{CalderoChapoton06} implies the identity
\begin{align}\label{eq:CR_counting_points}
\sum_{(\db,d):\mathrm{\db+\dd\ fixed}} \sharp(\Gr_\db F\mM)\sharp(\Gr_\dd F\mN) |\sK|^{ \langle  \dd,\dm-\db\rangle}= \sharp(\Gr_{\db+\dd} F\mM\oplus F\mN).
\end{align}
Now look at the desired identity \eqref{eq:spec_mult}. We have
\begin{align*}
LHS=&\sum_{\db,\dd} \sharp(\Gr_\db F\mM) |\sK|^{-\Hf  \langle \db,\dm-\db \rangle }X^{\ind_T \mM -\phi(\db)} \\
& \qquad \cdot \sharp(\Gr_\dd F\mN) |\sK|^{-\Hf  \langle \dd,\dn-\dd \rangle }X^{\ind_T\mN -\phi(\dd)}\\
=&\sum_{\db,\dd} \sharp(\Gr_\db F\mM)\sharp(\Gr_\dd F\mN) |\sK|^{-\Hf ( \langle \db,\dm-\db \rangle + \langle \dd,\dn-\dd \rangle )}\\
&\qquad \cdot |\sK|^{\Hf \lambda(\ind_T \mM -\phi(b),\ind_T \mN  -\phi(\dd))}X^{\ind_T \mM -\phi(b)+\ind_T \mN -\phi(\dd)}\\
=&|\sK|^{\Hf \lambda(\ind_T \mM ,\ind_T \mN )} \sum_{\db,\dd} \sharp(\Gr_\db F\mM)\sharp(\Gr_\dd F\mN) |\sK|^{  \langle  \dd,\dm-\db\rangle}\\
&\qquad \cdot |\sK|^{-\Hf  \langle \db+\dd,\dm+\dn-\db-\dd \rangle } X^{\ind_T \mM -\phi(b)+\ind_T \mN -\phi(\dd)}\\
=&|\sK|^{\Hf \lambda(\ind_T \mM ,\ind_T \mN )} \sum_{\db+\dd} \sharp(\Gr_{\db+\dd} F\mM\oplus F\mN) |\sK|^{-\Hf  \langle \db+\dd,\dm+\dn-\db-\dd \rangle } \\
&\qquad \cdot X^{\ind_T{ \mM \oplus \mN }-\phi(\db+\dd)}\\
=&RHS,
\end{align*}
where we have used Corollary \ref{cor:expand_lambda} for the third equality.

(ii) Assume that the objects $\mM$ and $\mN$ belong to $\add T$. The assertion follows from the definition of $\cT^\sK$.

(iii) Assume that the object $\mM$ contains no summands in $\add T$ and $\mN$ belongs to $\add T$. Then we have $\dd$=0, $\dn$=0. Write $\mN=\oplus_{1\leq i\leq n}T_i^{\dn_i}$, $\dn_i\geq 0$. The rigidity of $\mM\oplus \mN$ implies $\Hom_\cC(\mM,\tau \mN)=\Ext^1_\cC(\mM,\mN)=0$, and consequently the space $\Hom(F\mM,\oplus_{1\leq i\leq n} I_i^{\dn_i})$ vanishes. Then for any submodule $B$ of $F\mM$ with dimension vector $\db$, we have $\ind_T \mN\cdot \db=\sum_{1\leq i\leq n}\dn_i \db_i=0$. Starting from the left hand side of \eqref{eq:spec_mult} we obtain
\begin{align*}
LHS=&\sum_\db \sharp(\Gr_\db F\mM) |\sK|^{-\Hf  \langle \db,\dm-\db \rangle }X^{\ind_T \mM -\phi(\db)} X^{\ind_T\mN }\\
=&\sum_\db \sharp(\Gr_\db F\mM) |\sK|^{-\Hf  \langle \db,\dm-\db \rangle }\cdot |\sK|^{\Hf \lambda(\ind_T \mM -\phi(\db),\ind_T \mN )}\\
&\qquad X^{\ind_T \mM -\phi(\db)+\ind_T \mN }\\
=&|\sK|^{\Hf \lambda(\ind_T \mM ,\ind_T \mN )} \sum_\db  \sharp(\Gr_\db F\mM) |\sK|^{-\Hf  \langle \db,\dm-\db \rangle+\Hf \ind_T \mN \cdot \db }\\
&\qquad X^{\ind_T \mM -\phi(\db)+\ind_T \mN -\phi(\dd)}\\
=&|\sK|^{\Hf \lambda(\ind_T \mM ,\ind_T \mN )} \sum_\db  \sharp(\Gr_\db F(\mM\oplus \mN)) |\sK|^{-\Hf  \langle \db,\dm-\db \rangle }\\
&\qquad  X^{\ind_T{ \mM \oplus \mN }-\phi(\db+\dd)}\\
=&RHS.
\end{align*}

(iv) Assume that the object $\mN$ contains no summands in $\add T$ and that $\mM$ belongs to $\add T$. The proof is similar to case (iii).

(v) For any objects $\mM$, $\mN$ as in the proposition, write $\mM=\mM_1\oplus \mM_2$, $\mN=\mN_1\oplus \mN_2$, where $\mM_1,\mN_1$ contain no summands in $\add T$, and $\mM_2, \mN_2$ belong to $\add T$. Using the previous cases we obtain
\begin{align*}
X^\sK_\mM X^\sK_\mN=&|\sK|^{-\Hf \lambda([\mM_1],[\mM_2])-\Hf \lambda([\mN_1],[\mN_2])}X^\sK_{\mM_1} X^\sK_{\mM_2} X^\sK_{\mN_1} X^\sK_{\mN_2}\\
=&|\sK|^{-\Hf \lambda([\mM_1],[\mM_2])-\Hf \lambda([\mN_1],[\mN_2])}X^\sK_{\mM_1} |\sK|^{\lambda([\mM_2],[\mN_1])}X^\sK_{\mN_1} X^\sK_{\mM_2} X^\sK_{\mN_2}\\
=&|\sK|^{-\Hf \lambda([\mM_1],[\mM_2])-\Hf \lambda([\mN_1],[\mN_2])}|\sK|^{\lambda([\mM_2],[\mN_1])}\\
&\qquad \cdot |\sK|^{\Hf \lambda([\mM_1],[\mN_1])+\Hf \lambda([\mM_2],[\mN_2])}X^\sK_{\mM_1\oplus \mN_1} X^\sK_{\mM_2\oplus \mN_2}\\ 
=&|\sK|^{-\Hf \lambda([\mM_1],[\mM_2])-\Hf \lambda([\mN_1],[\mN_2])+\lambda([\mM_2],[\mN_1])+\Hf \lambda([\mM_1],[\mN_1])+\Hf \lambda([\mM_2],[\mN_2])}\\
&\qquad \cdot |\sK|^{\Hf \lambda([\mM_1]+[\mN_1],[\mM_2]+[\mN_2])}X^\sK_{\mM_1\oplus \mN_1\oplus \mM_2\oplus \mN_2}\\
=&|\sK|^{\Hf \lambda([\mM_1]+[\mM_2],[\mN_1]+[\mN_2])}X^\sK_{\mM\oplus \mN}\\
=&RHS.
\end{align*}
\end{proof}	
	
\begin{Rem}\label{rem:counting_polynomial_direct_sum}
Let $M$ and $L$ be given as in the above proposition. If for any submodule Grassmannians $\Gr_\db F\mM$ and $\Gr_\dd F\mN$, there exist some Laurent polynomial $P_\db$ and $Q_\dd$ in $q^\Hf$, \st over every finite field $\sK$ we have
\[
P_\db(|\sK|^\Hf)=\sharp(\Gr_\db F\mM),
\]
\[
Q_\dd(|\sK|^\Hf)=\sharp(\Gr_\dd F\mN),
\]
then it follows from equation \eqref{eq:CR_counting_points} that over every finite field $\sK$ we also have
\[
(\sum_{(\db,\dd):\db+\dd\, fixed}P_\db Q_\dd q^{  \langle  \dd,\dm-\db\rangle})|_{q^\Hf\mapsto |\sK|^\Hf}= \sharp(\Gr_{\db+\dd} F(\mM\oplus \mN)).
\]
\end{Rem}
By comparing this remark with Lemma \ref{lem:q_counting_polynomial} and Corollary \ref{cor:Katz_thm}, we obtain the following result.
\begin{Cor}\label{cor:Serre_polynomial_direct_sum}
If all submodule Grassmannians of $F\mM$ and $F\mN$ have their Serre polynomials as counting polynomials, then so do those of $F(\mM\oplus \mN)$. Moreover, for the elements $X_ \mM$, $X_ \mN$, $X_{ \mM \oplus  \mN }$ as in Definition \ref{def:X_M}, we have
\[
X_ \mM \cdot X_ \mN=q^{\Hf \lambda(\ind_T \mM ,\ind_T \mN )}X_{ \mM \oplus  \mN }.
\]
\end{Cor}	
	
\subsection{Specialized main theorem}
	
For convenience we fix the following conventions in this subsection:
\begin{itemize}
\item We assume that the ice quiver $\tQ$ is acyclic.
\item We work over a fixed finite field $\sK$. The cluster category $\cC$ is that of $\sK Q$, \cf Section \ref{sec:field_reduction}.
\item We denote an object of $\cC$ using a letter with a tilde, \eg $\tM\in \cC$, and denote its image under the functor 
\[
F=\Ext^1_\cC(T,\ )
\] by the same letter without a tilde, \eg $M=F(\tM)\in \mod \sK Q$. And we denote a morphism of $\cC$ and its image under $F$ by the same letter.
\item For any rigid object $\tL$ of $\cC=\cC_{\sK Q}$, we define $X^\sK_\tL$ and $\ind_T \tL$ to be the element in Definition \ref{def:spec_X_M} and the index associated with the lift of $\tL$ to $\cC_{k Q}$ respectively, \cf Part a) of Theorem \ref{thm:field_reduction}.
\item Recall that for $\sK Q$-modules $\mM$, $\mA$ and $\mB$, the \emph{Hall number} $F^\mM_{\mA\mB}$ is the number of submodules $\mU$ of $\mM$ \st $\mU$ is isomorphic to $\mB$ and $\mM/\mU$ is isomorphic to $\mA$.
\end{itemize}

	\begin{Prop}[Specialized mutation rule]\label{prop:specialized_mutation_rule}
Let $\tmM$ and $\tmN$ be two indecomposable coefficient-free rigid objects in $\cC$ \st the dimension of the space $\dim\Ext^1_\cC(\tmM,\tmN)$ equals $1$, and let
\begin{align}
&\tmN\ra \tmBp\ra \tmM\ra \tmN[1],\label{eq:1st_tri}\\
&\tmM\ra \tmBn\ra \tmN\ra \tmM[1]\label{eq:2nd_tri}
\end{align}
be two non-split triangles of $\cC$. If we have $\ind_T \tmN+\ind_T \tmM=\ind_T \tmBp$ as in Lemma \ref{lem:compute_lambda}, then the following identity holds:
\begin{align}\label{eq:desired_identity}
X^\sK_{\tmN}\cdot X^\sK_{\tmM}= |\sK|^{\Hf \lambda({\ind_T{\tmN}},{\ind_T{\tmBp}})}X^\sK_{\tmBp}+|\sK|^{\Hf \lambda({\ind_T{\tmN}},{\ind_T{\tmBn}})}X^\sK_{\tmBn}.
\end{align}
\end{Prop}
We postpone the proof for the moment. Notice that the (exchange) triangles \eqref{eq:1st_tri} and \eqref{eq:2nd_tri} of $\cC=\cC_{\sK Q}$ are compatible with those of $\cC_{kQ}$ by Part b) of Theorem \ref{thm:field_reduction}.

Recall that the ice quiver $\tQ$ is acyclic.
\begin{Thm}[Specialized main theorem]\label{thm:spec_X_M}
For any vertex $t$ of $\sT_n$ and any $1\leq i\leq n$, we have
\[
X^\sK_{T_i(t)}=X^\sK_i(t).
\]
Moreover, the map taking an object $M$ to $X^\sK_M$ induces a bijection from the set of isomorphism classes of coefficient-free rigid objects of $\cD$ to the set of quantum cluster monomials of $\cA^\sK$.
\end{Thm}
\begin{proof}
Recursively applying Proposition \ref{prop:specialized_mutation_rule} proves the assertion on quantum cluster variables. Next, inductively applying Proposition \ref{prop:direct_sum}, we prove the assertion on quantum cluster monomials.
\end{proof}
	
We shall divide Proposition \ref{prop:specialized_mutation_rule} into two cases and prove them in the rest of this subsection, using standard methods in representation theory and cluster theory as in \cite{Hubery} and \cite{CalderoKeller06}. We fix the following convention:
\begin{itemize}
\item We use the same symbol to denote a module and its isomorphism class.
\item As before, we shall use uppercase letters to denote $\sK Q$-modules and lowercase letters to denote their dimension vectors. We choose the the set formed by the $[S_i]$, $1\leq i\leq n$ as a basis of $\Kz(\mod \sK Q)$ and identify the dimension vectors with classes in $\Kz(\mod \sK Q)$.
\end{itemize}
		
	\subsubsection*{Proof in the generic cases}
		
		The following construction follows \cite[Section 7]{Hubery}. Assume that the objects $\tmM$ and $\tmN$ contain no summands of ${T}$. We write the indecomposable $\mM= \mM' \oplus \mP$, \st $ \mM' $ contains no projective summand, and $ \mP $ is some projective module. 
		
		Applying the functor $F=\Ext^1_\cC(T,\quad)$ to the triangle \eqref{eq:1st_tri}, we obtain a short exact sequence
\begin{align}\label{seq:generic_1st}
0\ra \mN \ra  \mBp  \ra \mM \ra 0.
\end{align}
For any submodule $H$ of $\mBp$, we have the following commutative diagram with exact rows and columns
\begin{align}\label{diag:1st_submod}
\begin{CD}
  @.		0     @.		0							@.				0\\
@.	    @AAA      	@AAA        						@AAA\\
0  @>>>     \mC     @>>>     \mG 			@>>>       \mA   @>>>  0\\
@.      @AAA        @AAA                 @AAA\\
0 @>>>   \mN 		@>>> 	 \mBp 	@>>> 			 \mM  @>>>  	0\\
@.	    @AAA 			  @AAA					@AAA\\
0 @>>>   \mD    @>>> 	  \mH   	   			@>>> 			\mB 	@>>> 	0.\\
@.			@AAA      	@AAA        						@AAA\\
  @.		0     @.		0							@.				0
\end{CD}
\end{align}
		
		Applying the functor $F=\Ext^1_\cC(T,\quad)$ to the triangle \eqref{eq:2nd_tri}, we obtain a long exact sequence in the module category
\begin{align}\label{seq:generic_2nd}
 \mM \xra{i}  \mW \oplus \mV  \xra{p} \mN \xra{\theta} \tau \mM ',
\end{align}
where $\mW$ and $\mV$ are defined by
\begin{align*}
& \mU =\Im\theta,\\
&0\ra  \mV \ra \mN \xra{\theta}  \mU \ra 0,\\
&0\ra  \mU \ra \tau  \mM' \ra\tau \mW '\oplus \mI\ra 0,\\
&\mW=\mW' \oplus \mP.
\end{align*}
Here $\mI$ is some injective module. One can prove
\begin{align*}
 \tmBn \cong { \tmW \oplus \tmV \oplus \tmI}[-1].
\end{align*}

Let $\mE$ be any submodule of $ \mW \oplus \mV $ and $\mF$ the quotient module $ (\mW \oplus \mV) /\mE$. We write $\mE=\mE_{\mL,\mY}$ and $\mF=\mF_{\mK,\mX}$, where the modules are defined by the following commutative diagram with exact rows and columns
\begin{align*}
\begin{split}
\begin{CD}
  @.		0     @.		0							@.				0\\
@.	    @AAA      	@AAA        						@AAA\\
0  @>>>     \mX     @>>>     \mF_{\mK,\mX} 			@>>>       \mK   @>>>  0\\
@.      @AAA        @AAA                 @AAA\\
0 @>>>   \mV 		@>>> 	 \mW \oplus  \mV 	@>>> 			 \mW  @>>>  	0\\
@.	    @AAA 			  @AAA					@AAA\\
0 @>>>   \mY    @>>> 	  \mE_{\mL,\mY}   	   			@>>> 			\mL 	@>>> 	0.\\
@.			@AAA      	@AAA        						@AAA\\
  @.		0     @.		0							@.				0
\end{CD}
\end{split}
\end{align*}

Now following \cite{CalderoKeller06} we can complete the long exact sequence \eqref{seq:generic_2nd} to a commutative diagram with exact rows and columns
\begin{align}\label{diag:submod}
\xymatrix{
		&  				&		0     		&		0								&		0				& 						&\\
		&	0\ar[r] 		&  		\mA\ar[r]\ar[u]    	&     \mF_{\mK,\mX}\ar[r]\ar[u]	&    \mC\ar[r]\ar[u]   	& \Cok p\ar[r] 			& 0\\
						&			 		&  \mM\ar[r]^i\ar[u]				&  \mW \oplus  \mV\ar[r]^p\ar[u] 	& \mN\ar[r]^{\theta}\ar[u] 			&  	\tau\mM'  &\\
0\ar[r]	& \Ker i\ar[r] 	&   \mB\ar[r]\ar[u]    		& 	  \mE_{\mL,\mY}\ar[r]\ar[u] & 	\mD\ar[r]\ar[u] 		& 	0, 				&\\
		&  				&		0\ar[u]     		&		0\ar[u]								&				0\ar[u]		& 							&
	}
\end{align}

where $\mB=i^{-1} (\mE_{\mL,\mY})$, $\mD=p(\mE_{\mL,\mY})$.
Notice that we have
\begin{align*}
\dim\Cok p&=\dim \mN-\dim \ker\theta=\dim \mN-\dim \mV=\dim \mU \\
\dim\Ker i&=\dim \mM+\dim \mN -(\dim  \mW\oplus \mV +\dim \Cok p)\\
&=\dim \mM-\dim  \mW =\dim  \mM' -\dim  \mW '
\end{align*}
Now look at the desired identity \eqref{eq:desired_identity}
\[
X^\sK_{\tmN}\cdot X^\sK_{\tmM}= |\sK|^{\Hf \lambda({\ind_T{\tmN}},{\ind_T{\tmBp}})}X^\sK_{\tmBp}+|\sK|^{\Hf \lambda({\ind_T{\tmN}},{\ind_T{\tmBn}})}X^\sK_{\tmBn}.
\]
Applying Corollary \ref{cor:expand_lambda}, we compute the left hand side as
\begin{align*}
LHS=&\sum_\dd \sharp(\Gr_\dd \mN) |\sK|^{-\Hf  \langle \dd,\dn-\dd \rangle } X^{{\ind_T \tmN}-\phi(\dd)}\\ 
&\qquad \sum_\db \sharp(\Gr_\db \mM) |\sK|^{-\Hf  \langle \db,\dm-\db \rangle } X^{{\ind_T \tmM}-\phi(\db)}\\
=&\sum_{\dd,\db} \sharp(\Gr_\dd \mN)\sharp(\Gr_\db \mM) |\sK|^{-\Hf ( \langle \dd,\dn-\dd \rangle + \langle \db,\dm-\db \rangle )} \\
&\qquad |\sK|^{\Hf \lambda({\ind_T \tmN}-\phi(\dd),{\ind_T \tmM}-\phi(\db))} X^{{\ind_T \tmN}-\phi(\dd)+{\ind_T \tmM}-\phi(\db)} \\
=&\sum_{\dd,\db} \sharp(\Gr_\dd \mN) \sharp(\Gr_\db \mM) |\sK|^{-\Hf ( \langle \dd,\dn-\dd \rangle + \langle \db,\dm-\db \rangle )}
|\sK|^{\Hf  (\langle \db,\dn \rangle -\langle \dd,\dm \rangle +\langle \dd,\db \rangle - \langle \db,\dd \rangle) }\\&\qquad |\sK|^{\Hf \lambda({\ind_T \tmN},{\ind_T \tmM})}X^{{\ind_T \tmN}-\phi(\dd)+{\ind_T \tmM}-\phi(\db)} \\
=&\sum_{\dd,\db} \sharp(\Gr_\dd \mN)\sharp(\Gr_\db \mM) |\sK|^{-\Hf  \langle \dd+\db,\dn+\dm-\dd-\db \rangle } |\sK|^{ \langle \db,\dn-\dd \rangle }\\
&\qquad |\sK|^{\Hf \lambda({\ind_T \tmN},{\ind_T \tmM})}X^{{\ind_T \tmN}-\phi(\dd)+{\ind_T \tmM}-\phi(\db)} .\\
\end{align*}
Notice that the number of rational points can be computed as $\sharp(\Gr_\dd \mN)=\sum_\mD F^\mN_{\mC\mD}$, $\sharp(\Gr_\db \mM)=\sum_\mB F^\mM_{\mA\mB}$. So we can write
\begin{align*}
LHS=&\sum_{\mB,\mD} F^\mN_{\mC\mD} F^\mM_{\mA\mB} |\sK|^{-\Hf  \langle \dd+\db,\dn+\dm-\dd-\db \rangle } |\sK|^{ \langle \mB,\mC \rangle }\\
&\qquad |\sK|^{\Hf \lambda({\ind_T \tmN},{\ind_T \tmM})}X^{{\ind_T \tmN}-\phi(\dd)+{\ind_T \tmM}-\phi(\db)}.
\end{align*}

Lemma \ref{lem:compute_lambda} implies that the first term of the right hand side is 
\begin{align*}
\Sigma_1=\sum_H F^{ \mBp }_{\mG\mH} |\sK|^{-\Hf  \langle \dh,\dg \rangle }|\sK|^{\Hf \lambda({\ind_T \tmN},{\ind_T \tmM})}X^{{\ind_T { \tmBp  }}-\phi(\dh)},
\end{align*}
and, since we have
\begin{align*}
\sharp(\Gr_{\dl+\dy} \mW \oplus \mV )=\sharp(\Gr_\dy  \mV ) \sharp(\Gr_\dl  \mW ) |\sK|^{ \langle \dl,\dv-\dy \rangle },
\end{align*}
we obtain that the 2nd term is
\begin{align*}
\Sigma_2=&\sum_{\mL,\mY} F^{ \mW }_{\mK\mL} F^{ \mV }_{\mX\mY} |\sK|^{ \langle \mL,\mX \rangle } |\sK|^{-\Hf  \langle \dl+\dy,\dw+\dv-\dl-\dy \rangle } |\sK|^{-\Hf }\\
&\qquad\cdot |\sK|^{\Hf \lambda({\ind_T {\tmN}},{\ind_T {\tmM}})}X^{{\ind_T { \tmBn }}-\phi(\dl)-\phi(\dy)}
\end{align*}

\begin{Lem}\label{lem:deg} For objects fitting into diagrams \eqref{diag:1st_submod} and \eqref{diag:submod}, we have
\begin{align*}
-\phi(\dh) + {\ind_T { \tmBp  }}=-\phi(\db+\dd) + {\ind_T \tmN} +{\ind_T \tmM}\\
-\phi(\dl+\dy) + {\ind_T { \tmBn }}=-\phi(\db+\dd) + {\ind_T \tmN} +{\ind_T \tmM}
\end{align*}
\end{Lem}
\begin{proof}
This is a well known result, \cf \cite[Lemma 16]{Palu08a}, \cite[Proposition 4.3]{CalderoKeller06}.
\end{proof}

According to Lemma \ref{lem:deg}, the desired identity \eqref{eq:desired_identity} follows from identities on Hall numbers, each of which takes the form
\begin{align}\label{eq:Hall_number_identity}
\begin{split}
&\sum_{\mB,\mD} F^\mN_{\mC\mD} F^\mM_{\mA\mB} |\sK|^{-\Hf  \langle \dd+\db,\dn+\dm-\dd-\db \rangle } |\sK|^{ \langle \mB,\mC \rangle } |\sK|^{\Hf \lambda({\ind_T \tmN},{\ind_T \tmM})}\\
=&\sum_H F^{ \mBp }_{\mG\mH} |\sK|^{-\Hf  \langle \dh,\dg \rangle }|\sK|^{\Hf \lambda({\ind_T \tmN},{\ind_T \tmM})} \\
&\qquad+ \sum_{\mL,\mY} F^{ \mW }_{\mK\mL} F^{ \mV }_{\mX\mY} |\sK|^{ \langle \mL,\mX \rangle } |\sK|^{-\Hf  \langle \dl+\dy,\dw+\dv-\dl-\dy \rangle } |\sK|^{-\Hf }|\sK|^{\Hf \lambda({\ind_T {\tmN}},{\ind_T {\tmM}})},
\end{split}
\end{align}
where all objects should fit into the same diagrams \eqref{diag:1st_submod} and \eqref{diag:submod}.

A similar identity on Hall numbers has been proved in \cite{Hubery} by using remarkable properties of Hall numbers, namely Green's Theorem \cite{Green95}\footnote{There is a gap in \cite{Hubery}: The existence of counting polynomials of Hall numbers is unknown. Nevertheless, the gap does not affect the identity on Hall numbers for a given finite field.}. We translate that identity to our setting as follows:
\begin{align}\label{eq:Hubery_identity}
\sum_{\mB,\mD} |\sK|^{ \langle \mB,\mC \rangle }F_{\mC\mD}^\mN F_{\mA\mB}^\mM = \sum_\mH F^{ \mBp }_{\mG\mH}+\sum_{\mL,\mY} |\sK|^{ \langle \dm-\dk,\dn-\dy \rangle }F^{ \mW }_{\mK\mL}F^{ \mV }_{\mX\mY}.
\end{align}
To show that \eqref{eq:Hall_number_identity} and \eqref{eq:Hubery_identity} are equivalent, it suffices to prove the following lemma\footnote{It might be possible to give an alternative proof of Proposition \ref{prop:specialized_mutation_rule}, by proving the following identity\begin{align*}
- \langle \dd,\da \rangle + \langle \db,\dc \rangle =&\Hf  \langle \db+\dd,\dm+\dn-\db-\dd \rangle \\
&\qquad -\Hf  \langle \dl+\dy,\dw+\dv-\dl-\dy \rangle -\Hf 
\end{align*}
We could then use the dichotomy phenomenon in the multiplication formula of cluster characters from \cite{Palu08a}. If we could follow this approach we would not need the result from \cite{Green95}.
}.
\begin{Lem}
One has the following identity
\begin{align}\label{eq:Euler_identity}
\begin{split}
- \langle \dl,\dx \rangle + \langle \dm-\dk,\dn-\dy \rangle =&\Hf  \langle \db+\dd,\dm+\dn-\db-\dd \rangle \\
&-\Hf  \langle \dl+\dy,\dw+\dv-\dl-\dy \rangle -\Hf .
\end{split}
\end{align}
\end{Lem}
\begin{proof}
Notice that both sides of the identity \eqref{eq:Euler_identity} are quadratic polynomials in the coordinates of $\de=\dl+\dy$ and $\dl$. Thus, we can expand the polynomials and compare the coefficients to prove the lemma.

First, expand both sides as
\begin{align*}
LHS=&- \langle  \mL, \mV -\de+\mL\rangle+ \langle \mM- \mW +\mL,\mN-\de+\mL \rangle \\
=& \langle \mL,\mN- \mV  \rangle + \langle \mM- \mW ,\mL \rangle - \langle \mM- \mW ,\de \rangle + \langle \mM- \mW ,\mN \rangle,
\end{align*}
\begin{align*}
RHS=&\Hf  \langle \mB+\mD,\mA+\mC \rangle -\Hf  \langle \de, \mW + \mV -\de \rangle -\Hf \\
=&\Hf  \langle \de+\Ker i, \mW + \mV +\Cok p-\de \rangle -\Hf  \langle \de, \mW + \mV -\de \rangle -\Hf \\
=&\Hf  \langle \Ker i, \mW + \mV +\Cok p \rangle +\Hf  \langle \de,\Cok p \rangle -\Hf  \langle \Ker i, \de \rangle -\Hf .\\
\end{align*}
Notice that we have
\begin{align*}
\dim \Ker i&=\dim \mM-\dim  \mW =\dim  \mM' -\dim  \mW ',\\
\dim \Cok p&=\dim \mN-\dim  \mV =\dim \mU,\\
\tau \dim \Ker i&= \dim \tau  \mM'  -\dim \tau \mW'  = \dim \mI + \dim \Cok p.
\end{align*}
In addition, for any dimension vector $x$, its Auslander-Reiten translation has the property
\[
\langle x,\  \rangle =- \langle \ ,\tau x \rangle.
\]
We can deduce the lemma from the following identities:

(i) $ \langle \mL,\mN- \mV  \rangle + \langle \mM- \mW ,\mL \rangle =0$.

Notice that we have
\[
\langle \mM- \mW ,\mL \rangle =- \langle \mL,\tau  \mM' -\tau { \mW '} \rangle =- \langle \mL, \mU +\mI \rangle.
\]
It suffices to show 
\[
\langle \mL,\mN- \mV - \mU  \rangle - \langle \mL,\mI \rangle =- \langle \mL,\mI \rangle =0.
\]
Since the object $\tmW \oplus \tmV \oplus \tmI[-1]$ is rigid, the Euler form 
$
 \langle  \mW ,\mI \rangle
$
vanishes,
and as a consequence we obtain
\[
\langle \mL,\mI \rangle =0.
\]

(ii) $- \langle \mM- \mW ,\de \rangle =\Hf  \langle \de,\Cok p \rangle -\Hf  \langle \Ker i, \de \rangle $.

Since $\Ker i=\mM- \mW$, it suffices to show 
\[
-\Hf  \langle \mM- \mW ,\de \rangle =\Hf  \langle \de,\Cok p \rangle.
\]
First, we have
\[
- \langle \mM- \mW ,\de \rangle = \langle \de,\tau  \mM' -\tau { \mW '} \rangle = \langle \de, \mU +\mI \rangle.
\]
Since the object $\tmW \oplus \tmV \oplus \tmI[-1]$ is rigid, the Euler form $\langle  \mW \oplus \mV ,\mI \rangle$ vanishes. Therefore we obtain 
\[
 \langle \de,\mI \rangle =0.
 \]
The identity then follows from $\dim\Cok p=\dim \mU $.

(iii) $ \langle \mM- \mW ,\mN \rangle =-1$.

Notice that we have
\begin{align*}
\langle \mM,\mN \rangle =&\dim \Hom(\mM,\mN)-\dim \Ext^1(\mM,\mN)\\
=&\dim \Hom(\mM,\mN)-1\\
=&\dim\Hom( \mM' ,\mN)+\dim\Hom(\mP,\mN)-1.
\end{align*}
From the rigidity of $\tmN\oplus\tmW \oplus \tmV \oplus \tmI[-1]$ we deduce
\[
 \langle  \mW ,\mN \rangle =\dim\Hom( \mW ,\mN)=\dim\Hom( \mW ',\mN)+\dim\Hom(\mP,\mN).
 \]
Consider the following exact sequences
\[
\Ext^1(\mN,\mN)\rightarrow \Ext^1(\mN, \mU )\rightarrow 0,
\]
\[
\Ext^1(\mN, \mU )\rightarrow \Ext^1(\mN,\tau  \mM' )\rightarrow \Ext^1(\mN,\tau \mW '\oplus \mI)\rightarrow 0.
\]
Notice that we have $\Ext^1(\mN,\tau  \mM' )\cong D\Hom( \mM' ,\mN)$, 
\[
\Ext^1(\mN,\tau  \mW ' \oplus \mI)=\Ext^1(\mN,\tau \mW ' )\cong D\Hom( \mW ',\mN),
\]
and also $\Ext^1(\mN,\mN)=0$. We obtain
\[
\Ext^1(\mN,\mU)=0,\]
\[
\dim \Hom( \mW ',\mN)=\dim \Hom( \mM' ,\mN).
\]
The desired identity follows.

(iv) $\langle \Ker i, \mW + \mV +\Cok p \rangle =-1$.

Apply $\tau$ to $\Ker i$. It suffices to show
\[
\langle \mN+ \mW , \mU +\mI \rangle =1.
\]
From the rigidity of $\tmW \oplus \tmV \oplus \tmI[-1]\oplus \tmN$, we deduce 
\[
\langle \mN+ \mW ,\mI \rangle =0,
\]
\[
\langle \mN+ \mW , \mU  \rangle = \langle \mN, \mU  \rangle + \langle  \mW , \mU  \rangle.
\]
The identity (iv) follows from the following two assertions.

(iv-a) $ \langle \mN, \mU  \rangle =1$.

Notice that we have a sequence
\[
0\rightarrow \Hom (\mN, \mU )\rightarrow \Hom(\mN,\tau  \mM' ),
\]
and also $\dim \Hom(\mN,\tau  \mM' )=1$, $0\neq \theta\in \Hom(\mN, \mU )$. We deduce
\[
\dim \Hom(\mN, \mU )=1.
\]
Consider the sequence
\[
\Ext^1(\mN,\mN)\rightarrow \Ext^1(\mN, \mU )\rightarrow 0,
\]
Since the object $\mN$ is rigid, the space $\Ext^1(\mN, \mU )$ vanishes.

Therefore, we obtain 
\[
\langle \mN, \mU  \rangle =\dim\Hom (\mN, \mU )-\dim\Ext^1(\mN, \mU )=1.
\]

(iv-b) $ \langle  \mW , \mU  \rangle =0$.

Notice that we have a sequence
\[
\Ext^1( \mW ,\mN)\rightarrow \Ext^1( \mW , \mU )\rightarrow 0,
\]
and also the rigidity of $\tmW \oplus \tmV \oplus \tmI[-1]\oplus \tmN$. We deduce $\Ext^1( \mW ,\mN)=0$, and further $\Ext^1( \mW ,\mU)=0$.

Consider the sequence
\[
\Ext^1(\tau  \mM' ,\tau  \mW ')\rightarrow \Ext^1( \mU ,\tau \mW' )\rightarrow 0.
\]
Since the space $\Ext^1(\tau  \mM' ,\tau  \mW ')$ vanishes, we find $D\Hom( \mW ',\mU)\cong\Ext^1(\mU,\tau  \mW ')$ to be zero. Finally consider the sequence 
\[
0\ra \Hom(\mP, \mU )\ra \Hom(\mP,\tau  \mM' )\cong \mD \Ext^1( \mM' ,\mP)\subset\Ext^1(\mM,\mM)=0.
\]
We obtain
\[
\dim\Hom(\mW,\mU)=\dim\Hom(\mW' ,\mU)+\dim\Hom(\mP,\mU)=0.
\]
Therefore, we have
\[
\langle  \mW , \mU  \rangle=\dim\Hom(\mW,\mU)-\dim\Ext^1(\mW,\mU)=0.
\]

\end{proof}
	
\subsubsection*{Proof in the exceptional cases}

If both objects $\tmN$ and $\tmM$ belong to $\add T$, then $\Ext^1_\cC(\tmN,\tmM)\neq 0$ contradicts $\Ext^1_\cC(T,T)= 0$.

If $\tmM$ belongs to $\add T$ and $\tmN$ is not in $\add T$, by lifting the triangle \eqref{eq:1st_tri} to $\cF$ we obtain a non-split triangle
\[
\pi^{-1}\tmN\ra \pi^{-1}\tmBp\ra \pi^{-1}\tmM\ra \pi^{-1}\tmN[1].\\
\]
It follows that $\Hom_{\per\Gamma}(\pi^{-1}\tmM,\pi^{-1}\tmN[1])$ is nonzero. But since the object $\pi^{-1}\tmM$ belongs to $\add\Gamma$ and $\pi^{-1}\tmN$ belongs to $\pr_{\per\Gamma}\Gamma$, this contradicts the fact that $\Hom_{\per\Gamma}(\Gamma,\Gamma[i])$ vanishes for any $i>0$.

So we are left with the case when $\tmN=T_k$ for some $1\leq k\leq n$ and $\tmM\notin\add T$. There are two non-split triangles in the cluster category
\begin{align}
 T_k\rightarrow  \tmBp \rightarrow \tmM (\xrightarrow{\mathsf {\theta}} \tmI_k),\label{tri:spec_1st}\\
(\tmP_k\xrightarrow{\mathsf{\phi}})\tmM\rightarrow  \tmBn \rightarrow  T_k. \label{tri:spec_2nd}
\end{align}
In addition, Lemma \ref{lem:compute_lambda} implies the identities
\begin{align*}
\Hf \lambda(\ind_T T_k,{\ind_T { \tmBp }})=&\Hf \lambda(\ind_T T_k,{\ind_T {\tmM}}),\\
\Hf \lambda(\ind_T T_k,{\ind_T { \tmBn }})=&\Hf \lambda(\ind_T T_k,{\ind_T {\tmM}})-\Hf .
\end{align*}
By applying the functor $F$ to the triangles \eqref{tri:spec_1st} and \eqref{tri:spec_2nd}, we obtain two exact sequences
\begin{align}
0\rightarrow  \mBp \rightarrow \mM \xrightarrow{\mathsf {\theta}} \mI_k,\label{seq:spec_1st}\\
\mP_k\xrightarrow{\mathsf{\phi}}\mM\rightarrow  \mBn \rightarrow  0. \label{seq:spec_2nd}
\end{align}
Following \cite{Hubery}, define $\mV=\Ker\theta$, $\mI=\Cok\theta$, $\mU=\Cok\phi$, $\mP=\Ker\phi$, and then we have
\begin{align*}
 \tmBp &=\tmV\oplus \tmI[-1],\\
 \tmBn &=\tmU \oplus \tmP[1].
\end{align*}
%Set $\mTheta=\Im \theta$, $\mPhi=\Im\phi$.
Now look at the desired identity \eqref{eq:desired_identity}
\begin{align*}
X^\sK_k X^\sK_{\tmM} = |\sK|^{\Hf \lambda(\ind_T T_k,{\ind_T { \tmBp }})}X^\sK_{ \tmBp }+|\sK|^{\Hf \lambda(\ind_T T_k,{\ind_T { \tmBn }})}X^\sK_{ \tmBn }.
\end{align*}
We expand both sides as
\begin{align*}
LHS=&|\sK|^{\Hf \lambda(\ind_T T_k,{\ind_T {\tmM}})}\sum_\mB F^\mM_{\mA\mB}|\sK|^{-\Hf  \langle \mB,\mA \rangle } |\sK|^{\Hf \lambda(\ind_T T_k,-\phi(\db))}\\
&\qquad X^{{\ind_T {\tmM}}-\phi(\db)+\ind_T T_k}\\
=&|\sK|^{\Hf \lambda(\ind_T T_k,{\ind_T {\tmM}})}\sum_\mB F^\mM_{\mA\mB}|\sK|^{-\Hf  \langle \mB,\mA \rangle } |\sK|^{-\Hf b_k}X^{{\ind_T {\tmM}}-\phi(\db)+\ind_T T_k},
\end{align*}
\begin{align*}
RHS&=|\sK|^{\Hf \lambda(\ind_T T_k,{\ind_T {\tmM}})}X^\sK_{ \tmBp }+|\sK|^{\Hf \lambda(\ind_T T_k,{\ind_T {\tmM}})}|\sK|^{-\Hf }X^\sK_{ \tmBn }\\
&=|\sK|^{\Hf \lambda(\ind_T T_k,{\ind_T {\tmM}})}\sum_\mD F^{\mV}_{\mC\mD}|\sK|^{-\Hf  \langle \mD,\mC \rangle } X^{{\ind_T {\tmV}}+{\ind_T {\tmI[-1]}}-\phi(\dd)}\\
&+|\sK|^{\Hf \lambda(\ind_T T_k,{\ind_T {\tmM}})}|\sK|^{-\Hf }\sum_\mF F^{ \mU }_{\mE\mF}|\sK|^{-\Hf  \langle \mF,\mE \rangle } X^{{\ind_T \tmU }+{\ind_T {\tmP[1]}}-\phi(\df)}.
\end{align*}
According to \cite[Section 3.2, Lemma 1]{CalderoKeller06} or \cite{Palu08a}, there is a dichotomy phenomenon in the module category:
\begin{enumerate}
\item When $\db_k$ $(=(\dim \mB)_k)$ vanishes, we have
\[
\mD\cong \mB,\quad (\dim \mD)_k=0,
\]
\item When $\db_k$ equals $1$, we have
\[
0\rightarrow \mPhi\rightarrow \mB \rightarrow \mF\rightarrow 0.
\]
Moreover, since $\phi$ is non trivial, we have $(\dim \mPhi)_k=1$.

\end{enumerate}

It suffices to show that in case (1) we have the identity
\begin{align*}
-\Hf  \langle \mB,\mA \rangle -\Hf \db_k=-\Hf  \langle \mD,\mC \rangle ,
\end{align*}

and in case (2) the identity
\begin{align*}
-\Hf  \langle \mB,\mA \rangle -\Hf \db_k=-\Hf  \langle \mF,\mE \rangle -\Hf .
\end{align*}

\begin{itemize}
\item Proof in case (1):
From the facts
\begin{align*}
 \langle \mB,\mA \rangle - \langle \mD,\mC \rangle = &\langle \mD,\mM-\mD \rangle - \langle \mD,\mV-\mD \rangle \\
= &\langle \mD,\mM-\mV \rangle = \langle \mD,\mTheta \rangle,
\end{align*}
\[
0\rightarrow \Hom(\mD,\mTheta)\rightarrow \Hom(\mD,\mI_k),
\]
\[
\Hom(\mD,\mI_k)=\Hom(\mB,\mI_k)=0,
\]
we deduce that $ \langle \mB,\mA \rangle - \langle \mD,\mC \rangle$ vanishes.
\item Proof of case (2):
Notice that by applying the functor $F$ to the triangle \eqref{tri:spec_2nd} we get a commutative diagram
\[
\xymatrix{
         								& 0  														& 0   						&\\
         								& \mA \ar[u]\ar[r]							& \mE\ar[u]\ar[r]   &0\\
 \mP_k \ar[r]^{\phi}			& \mM  \ar[u]\ar[r]							& \mU\ar[u]\ar[r]  &0\\
 \mP_k\ar@{=}[u]\ar[r] 		& \mB  \ar[u]\ar[r]							& \mF\ar[u]\ar[r]   &0\\
 												&	0			\ar[u]									& 0. \ar[u]				&
 }
\]
Then the module $A$ is isomorphic to $N$. As a consequence, we have
\begin{align*}
 \langle \mB,\mA \rangle - \langle \mF,\mE \rangle &= \langle \mB-\mF,\mE\rangle\\
 &=\langle \mPhi,\mE\rangle.
\end{align*}
Consider the short exact sequence
\[
0\rightarrow \Hom(\mP_k,{\mPhi})\rightarrow\Hom(\mP_k,\mM)\rightarrow \Hom(\mP_k, \mU )\rightarrow 0.
\]
Since $\Hom(\mP_k,\mPhi)$ does not vanish and $\Hom(\mP_k,\mM)$ is a subset of $\Hom_\cC(\tau^{-1} T_k,\tmM)$ whose dimension is $1$, we obtain that $\Hom(\mP_k, \mU )$ vanishes.

Also, the rigidity of $T_k\oplus \tmBn$ implies that $\langle \mP_k, \mE\rangle$ vanishes. Consequently we have
\[
 \langle \mPhi,\mE \rangle = \langle \mP_k,\mE \rangle - \langle \mP,\mE \rangle =- \langle \mP,\mE \rangle.
 \]
The exact sequence
\[
0\rightarrow \Hom(\mP,\mF)\rightarrow \Hom(\mP, \mU )\rightarrow \Hom(\mP,\mE)\rightarrow 0
\]
and the inclusion of vector spaces
\[
\Hom(\mP, \mU )\subset \Hom_{\cC}(\tmP[1]\oplus \tmU,\Sigma (\tmP[1]\oplus \tmU))=0
\]
imply that $\langle \mP,\mE \rangle$ vanishes. Therefore we have
\begin{align*}
\langle \mB,\mA \rangle - \langle \mF,\mE \rangle&= \langle \mPhi,\mE \rangle \\
&=-\langle \mP,\mE \rangle\\
&=0.
 \end{align*}
\end{itemize}

\section{Proof of Quantum \texorpdfstring{$F$}{F}-polynomials}\label{sec:proof_F_M}

%\input{proof_F_M.tex}
%\section{Proof of Quantum \texorpdfstring{$F$}{F}-polynomials}\label{sec:proof_F_M}

In this section, we assume that the ice quiver $\tQ$ is acyclic and show how Theorem \ref{thm:F_M} follows from Theorem \ref{thm:acyclic_main_thm}.

Because $\tQ$ is acyclic, there is a bijection between the set of coefficient-free rigid objects in $\cC$ and rigid objects in $\cC_Q$, \cf Section 4 of \cite{IyamaYoshino08}. We denote coefficient-free rigid objects in $\cC$ and their images under the bijection by the same symbols. For any vertex $t$ of $\sT_n$ and any $1\leq i\leq n$, according to Theorem \ref{thm:acyclic_main_thm} and Proposition \ref{prop:decompose_X_M} we have
\begin{align*}
X_i(t)&=\sum_e E(\Gr_e \Ext^1_{\cC}(T,T_i(t))) q^{-\Hf \langle e,[\Ext^1_{\cC}(T,T_i(t))]-e \rangle} X^{\ind_T T_i(t)-\phi(e)}\\
&=\sum_e E(\Gr_e \Ext^1_{\cC}(T,T_i(t))) q^{\Hf \langle e,e \rangle}X^{-\phi(e)} X^{\ind_T T_i(t)}.
\end{align*}
Theorem \ref{thm:Tran_thm} yields that
\[
X_i(t)=F_i(t)|_{y^e\mapsto X^{-\phi(e)}}X^{\ind_T T_i(t)}.
\]
Therefore, the quantum $F$-polynomial $F_i(t)$ equals $F_{T_i(t)}$.

\section{Proof for arbitrary coefficients}\label{sec:proof_general_coefficient}

%\input{proof_general_coefficient.tex}
%\section{Proof for arbitrary coefficients}\label{sec:proof_general_coefficient}

This section is devoted to completing the proof of the main Theorem \ref{thm:X_M}. For any vertex $t$ in $\sT_n$ and any $1\leq i\leq n$, denote by $F_i^\Z(t)$ the commutative $F$-polynomial associated to $i$ and $t$, and by $F_i(t)$ the quantum one. According to Theorem \ref{thm:F_M}, the element $F_{T_i(t)}$ belongs to the set of quantum $F$-polynomials. We consider the following map
\[
\begin{matrix}
\{F_i(t),\ t\neq t_0\}&\xRa{f}&\{F_i^\Z(t),\ t\neq t_0\}\\
 F_{T_i(t)}&\longmapsto &ev_1(F_{T_i}(t))\\
 &&\parallel\\
 F_i(t)&\longmapsto&F_i^\Z(t).
\end{matrix}
\]
Because different nontrivial $F$-polynomials have different highest degrees, the map $f$ is bijective. Thus, we obtain
\[
F_{T_i(t)}=F_i(t)
\]
for $t\neq t_0$. Further, this identity holds for $t=t_0$ as well.

By Proposition \ref{prop:decompose_X_M}, the map in Theorem \ref{thm:X_M} takes $T_i(t)$ to the element
\[
X_{T_i(t)}=F_{T_i(t)}X^{\ind_T T_i(t)}.
\]
Thus, we have 
\[
X_{T_i(t)}=F_{T_i(t)}X^{\ind_T T_i(t)}=F_i(t)X^{\ind_T T_i(t)}=X_i(t)
\]
for any $1\leq i\leq n$ and any vertex $t$ of $\sT_n$, where the last equality comes from Theorem \ref{thm:Tran_thm}.

Next, if we inductively apply Corollary \ref{cor:Serre_polynomial_direct_sum} which claim that $X_?$ is multiplicative under taking direct sums, we obtain Theorem \ref{thm:X_M}.

%%%%%%%%%%%%%%%%%%%%%%%%%%
%%				Appendix
%%%%%%%%%%%%%%%%%%%%%%%%%%

%% The Appendices part is started with the command \appendix;

%% appendix sections are then done as normal sections
\appendix
\section{On the integral cluster category, \texorpdfstring{by Bernhard Keller}{by Bernhard Keller}}

\subsection{Integrality of rigid modules} \label{sec:int_rigid_mod}

Let $Q$ be a finite quiver
without oriented cycles. Recall that a {\em $\Z Q$-lattice} is a
$\Z Q$-module which is finitely generated free over $\Z$.
A $\Z Q$-lattice $M_0$ is {\em rigid} if the group of its selfextensions
vanishes. Part a) of the following theorem is a consequence
of Crawley-Boevey's results in \cite{CrawleyBoevey96};
part b) in addition uses Schofield's \cite{Schofield92},
cf. section~4 of \cite{CrawleyBoevey96}.

\begin{Thm}[Crawley-Boevey] \label{thm:class-rigid-lattices}

a) The functor which takes $M_0$ to $M_0\ten_{\Z} \C$ induces
a bijection from the set of isomorphism classes of rigid
$\Z Q$-lattices to the set of isomorphism classes of rigid
$\C Q$-modules; indecomposables correspond to indecomposables.

b) For each finite field $\K$, the functor which takes $M_0$ to
$M_0\ten_{\Z} \K$ induces a bijection from the set of isomorphism
classes of rigid $\Z Q$-modules to the set of isomorphism classes
of rigid $\K Q$-modules; indecomposables correspond to indecomposables.
\end{Thm}

The theorem shows in particular that for each rigid $\C Q$-module $M$
and each dimension vector $e$, the Grassmannian of submodules
$Gr_e(M)$ is canonically the extension to $\C$ of a scheme 
defined over $\Z$, namely, if $M_0$ is a rigid $\Z Q$-lattice
with an isomorphism from $M_0 \ten_{\Z} \C$  to $M$, then it
is the scheme $Gr_e(M_0)$ whose $\Z$-points are the $\Z Q$-sublattices
$U_0$ of $M_0$ which have rank vector $e$ and whose inclusion
admits a $\Z$-linear retraction. 

% For a finite field $\K$, we will sometimes
% write $Gr_e(M)(\K)$ for the $\K$-points of the scheme $Gr_e(M_0)$ and,
% if $\K$ is clear from the context, we write $\# Gr_e(M)$ for the number of
% $\K$-points of $Gr_e(M_0)$.

\subsection{Extensions between lattices}
Let $\lat \Z Q$ denote the category of $\Z Q$-lattices. It becomes
an exact category when we endow it with the class of exact sequences
whose underlying sequences of abelian groups are split exact. 

\begin{Lem} \label{lemma:ext-lattices}
Let $L$ and $M$ be $\Z Q$-lattices and $\F$ a field.
\begin{itemize}
\item[a)] For each field $\F$, 
the canonical morphism in the derived category of $\F$-vector spaces
\begin{equation} \label{eq:red-isom}
\RHom_{\Z Q}(L,M) \lten_Z \F \to \RHom_{\F Q}(L\ten_\Z \F, M\ten_\Z \F)
\end{equation}
is invertible.
\item[b)] If $L$ and $M$ are rigid, the groups $\Hom_{\Z Q}(L,M)$ and
$Ext^1_{\Z Q}(L,M)$ are finitely generated free and the canonical
maps 
\begin{align*}
\Hom_{\Z Q}(L,M)\ten_\Z \F & \to \Hom_{\F Q}(L\ten_\Z \F, M\ten_\Z \F) \quad
\mbox{and}  \quad \\
\Ext^1_{\Z Q}(L,M)\ten_\Z \F  & \to \Ext^1_{\F Q}(L\ten_\Z \F, M\ten_\Z \F)
\end{align*}
are isomorphisms. 
\end{itemize}
\end{Lem}

\begin{proof}  
a) Let us call {\em stalk lattice} a lattice which is concentrated
at a single vertex of $Q$. 
Since $Q$ does not have oriented cycles, each lattice has
a (finite) filtration whose subquotients are stalk lattices
of rank one and whose inclusions admit $\Z$-linear retractions.  
Clearly the assertion holds for stalk lattices of rank one and is
preserved under taking extensions in the first or the second argument.
Thus, it holds for all lattices.

b) By Lemma~2 of \cite{CrawleyBoevey96}, the first assertion
holds. The remark immediately preceding Lemma~2 in [loc. cit.] 
also yields that the second canonical map is an isomorphism.
By part a), it follows that the first canonical map is an isomorphism.
\end{proof}

\subsection{Reduction of cluster categories} 
\label{sec:field_reduction}
Let $\cd^b(\lat \Z Q)$ denote
the bounded derived category of the exact category $\lat \Z Q$,
cf. \cite{Keller96}. It is a triangulated category. Let $S$ denote
the total derived functor of the tensor product with the bimodule
$\Hom_\Z(\Z Q, \Z)$. Then $S$ is a Serre functor for $\cd^b(\lat \Z Q)$
in the sense that we have canonical bifunctorial isomorphisms
\[
\RHom_\Z(\RHom_{\Z Q}(L,M),\Z) \iso \RHom_{\Z Q}(M, SL).
\]
Clearly, for each field $\F$ and for $L$ in $\cd^b(\lat \Z Q)$, we
have a functorial isomorphism
\begin{equation} \label{eq:red-Serre}
(SL)\ten_\Z \F \iso S(L \ten_\Z \F) \; ,
\end{equation}
where the $S$ on the right denotes the Serre functor of the
bounded derived category of $\F Q$. We define the {\em integral
cluster orbit category $\cc^{orb}_{\Z Q}$} to be the orbit category of 
$\cd^b(\lat \Z Q)$ under the action of the automorphism
$\Sigma^{-2} S$, where $\Sigma$ is the suspension functor,
cf. \cite{Keller05}. Clearly, the
suspension functor of $\cd^b(\lat \Z Q)$ induces an auto-equivalence
of $\cc^{orb}_{\Z Q}$. We still denote it by $\Sigma$. 
We do not know whether the integral cluster orbit category
is triangulated in general and we suspect that this is not
the case. However, we can embed it into a
triangulated category using the technique of \cite{Amiot09}: 
Let $A$ be the $3$-derived preprojective algebra $\Pi_3(\Z Q)$
as defined in \cite{Keller09a}. By theorem~6.3 of
\cite{Keller09a}, we can describe the algebra $A$ 
explicitly as the Ginzburg dg algebra $\Gamma(\Z Q, 0)$.
We define $\cd_{per(\Z)}(A)$ to be be full subcategory
of the derived category $\cd(A)$ whose objects are the
dg $A$-modules $M$ whose underlying complex of
$\Z$-modules is perfect. Since $A$ is homologically
smooth, the category $\cd_{\per(\Z)}(A)$ is a full
subcategory of $\per(A)$.  We define the {\em integral cluster
category $\cc_{\Z Q}$} as the triangle quotient
\[
\cc_{\Z Q} = \per(A)/\cd_{\per(\Z)}(A).
\]
\begin{Prop} \label{prop:orbit-into-tria}
The tensor product $?\ten_{\Z Q} A$ induces
a fully faithful functor
\[
\cc^{orb}_{\Z Q} \to \cc_{\Z Q}
\]
and the category $\cc_{\Z Q}$ is the closure of its
image under extensions.
\end{Prop}

\begin{proof} The proof uses the technique of Theorem~7.1 of 
\cite{Keller05}. Details will be given elsewhere.
\end{proof}

For two objects $L$ and $M$ of the integral cluster category, we
define
\[
\Ext_{\cc_{\Z Q}}^i(L,M) = \Hom(L, \Sigma^i M)
\]
for all integers $i$. 

The description of the classical cluster category via
the Ginzburg algebra given in \cite{Amiot09} shows
that  for each field $\F$, the
extension of scalars from $\Z$ to $\F$ induces a canonical
additive functor, the {\em reduction functor}
\[
?\ten_\Z \F: \cc_{\Z Q} \to \cc_{\F Q}.
\]
Let us call an object of $\cc_{\Z Q}$ {\em lattice-like}
if it is a direct sum of images of direct summands of $\Sigma \Z Q$
and of images of $\Z Q$-lattices.
\begin{Lem}
\begin{itemize}
\item[a)] For two rigid lattices $L$ and $M$, we have a canonical
isomorphism
\[
\Ext^1_{\cc_{\Z Q}}(L,M) \iso \Ext^1_{\Z Q}(L,M) \oplus D_\Z \Ext^1_{\Z Q}(M,L) \; ,
\]
where $D_\Z$ is the functor $\Hom_\Z(?,\Z)$. 
\item[b)] 
For two rigid lattice-like objects $L$ and $M$ of $\cc_{\Z Q}$, the
group $\Ext^1(L,M)$ is finitely generated free and for each
field $\F$, we have the isomorphism
\[
\Ext^1(L,M) \ten_\Z \F \iso \Ext^1(L\ten_\Z \F, M\ten_\Z \F).
\]
\end{itemize}
\end{Lem}

\begin{proof} a) Thanks to proposition~\ref{prop:orbit-into-tria},
we can use the orbit category to compute the
morphisms in the cluster category. This computation is similar to the classical
one in \cite{BuanMarshReinekeReitenTodorov06}.
Part b) follows from a) and from part b) of
Lemma~\ref{lemma:ext-lattices}.
\end{proof}

An object $L$ of $\cc_{\Z Q}$ is {\em rigid} if the
group $\Ext^1(L,L)=\Hom(L,\Sigma L)$ vanishes.
Two objects $L$ and $M$ of $\cc_{\Z Q}$ form an
{\em exchange pair} if they are indecomposable and
rigid and the space $\Ext^1(L,M)$ is of rank one.

\begin{Thm} \label{thm:field_reduction}
a)  The {\em reduction functor}
\[
?\ten_\Z \F: \cc_{\Z Q} \to \cc_{\F Q}.
\]
induces a bijection from
the set of isomorphism classes of lattice-like rigid objects in $\cc_{\Z Q}$ to 
the set of isomorphism classes of rigid objects of $\cc_{\F Q}$. Under this bijection,
indecomposables correspond to indecomposables.

b) Two indecomposable lattice-like rigid objects $L$ and $M$
of $\cc_{\Z Q}$ form an exchange pair if and only if this holds for their images
in $\cc_{\F Q}$. In this case, one of the groups $\Ext^1_{\Z Q}(L,M)$  and
$\Ext^1_{\Z Q}(M,L)$ is free of rank one and the other vanishes.
If $\Ext^1_{\Z Q}(M,L)$ is free of rank one, there are triangles
\[
L \to E \to M \to \Sigma L \mbox{ and }  M \to E' \to L \to \Sigma M
\]
in $\cc_{\Z Q}$ whose images under the reduction functor
\[
\cc_{\Z Q} \to \cc_{\F Q}
\]
are the exchange triangles associated with the exchange pair
formed by $L\ten_\Z \F$ and $M\ten_\Z \F$.
\end{Thm}

\begin{proof} Part a) follows from the classification
of rigid lattices in Theorem~\ref{thm:class-rigid-lattices} and
the computation of extension groups in Lemma~\ref{lemma:ext-lattices}.
Part~b) follows from a) and Lemma~\ref{lemma:ext-lattices}.
\end{proof}

\subsection*{Acknowledgment} The author thanks Sarah Scherotzke
for pointing out a mistake in a previous version of 
lemma~ \ref{lemma:ext-lattices}

%%%%%%%%%%%%%%%%%%%%%%%%%%
%%			Bibliography
%%%%%%%%%%%%%%%%%%%%%%%%%%

%% If you have bibdatabase file and want bibtex to generate the
%% bibitems, please use
%%
%%  \bibliographystyle{elsarticle-harv} 
%%  \bibliography{  \langle  your bibdatabase \rangle }
%\bibliographystyle{amsalpha}%Can also use halpha.bst.
%\bibliography{referenceURL}
\newcommand{\etalchar}[1]{$^{#1}$}
\def\cprime{$'$} \def\cprime{$'$}
\providecommand{\bysame}{\leavevmode\hbox to3em{\hrulefill}\thinspace}
\providecommand{\MR}{\relax\ifhmode\unskip\space\fi MR }
% \MRhref is called by the amsart/book/proc definition of \MR.
\providecommand{\MRhref}[2]{%
  \href{http://www.ams.org/mathscinet-getitem?mr=#1}{#2}
}
\providecommand{\href}[2]{#2}

%% else use the following coding to input the bibitems directly in the
%% TeX file.

%\begin{thebibliography}{00}

%% \bibitem[Author(year)]{label}
%% Text of bibliographic item

%\bibitem[ ()]{}

%\input{ref.tex}
%\bibitem[Kontsevich(2008)]{KontsevichSoibelman08}
%why doesn't this latex code work?

%\end{thebibliography}

%%%%%%%%%%%%%%%%%%%%%%%%%
%% List of modifications
%%%%%%%%%%%%%%%%%%%%%%%%%
%\input{modification.tex}

\end{document}